\documentstyle[11pt]{article}

\def\Box{\rule{.5em}{.5em}}

\input amssym.def
\input amssym.tex

\newtheorem{theorem}{{Theorem}}[section]
\newtheorem{proposition}[theorem]{{Proposition}} 
\newtheorem{definition}[theorem]{{Definition}} 
\newtheorem{lemma}[theorem]{{Lemma}} 
\newtheorem{corollary}[theorem]{{Corollary}} 
\newtheorem{fact}[theorem]{{Fact}} 
\newtheorem{remark}[theorem]{{\it Remark}} 
\newtheorem{question}[theorem]{{Question}} 

\newenvironment{proof}{{\it Proof.}}{\hfill$\Box$\medskip}

\title{Isometry groups and geodesic foliations of Lorentz
manifolds. \\ Part I: Foundations of  Lorentz dynamics }
\author{Abdelghani Zeghib}
\date{\today}

\begin{document}

\maketitle

\begin{abstract} This is the first part of a series on
 non-compact
 groups acting isometrically on compact Lorentz manifolds.
This subject was recently  investigated by many authors.
In
the present part we investigate the dynamics of affine, and
 especially  Lorentz transformations. In
particular
we show how this is related to
 geodesic foliations. The existence of geodesic
foliations
was (very succinctly) mentioned  for   the first time
 by D'Ambra and Gromov, who suggested that
this may help in the classification
of  compact Lorentz manifolds with non-compact isometry
groups.
 In the Part II of the series, a partial classification
of compact Lorentz manifolds with
non-compact isometry group will
be achieved
 with the aid  of
 geometrical tools along  with the dynamical ones presented
 here.

\end{abstract}



\section{Lorentz dynamics}

We are interested  here in    the following question:
 when is the
isometry group
of a compact Lorentz
manifold  non-compact?
Many progress has been made towards answering it, see
for instance the works by, Zimmer \cite{Zim},
Gromov (\cite{Gro}, \cite{D-G}),
 D'Ambra \cite{D'A}, Kowalsky \cite{Kow}, Adams-Stuck
(\cite{A-S.1}, \cite{A-S.2}), and
the author (\cite{Zeg.espacetemps.homogenes},
\cite{Zeg.identity.component}).
At this stage, from the different  investigations
due to these authors, we
know the list of non-compact
{\bf connected} Lie groups acting isometrically
on compact Lorentz manifolds.   For some  groups
in the list, we also understand completely  the geometric
structure of the underlying  Lorentz manifold.
 However, we do not know  enough things about this structure,
in the case of the remaining groups, such as,  for example,
 the case of abelian groups.
 In the case of  {\bf non-connected} groups,
 no significant thing is known on neither their
algebraic structure, nor the geometric structure
of the Lorentz manifold on which they act.

Our approach here, is to study the
dynamics of lorentz transformations, i.e. diffeomorphisms
preserving Lorentz metrics. By this approach, we intend
to understand,  at the same time,  the structure of Lorentz manifolds
with large isometry groups, and the dynamics of the
individual isometries themselves.

It is natural to think about Lorentz metrics as the nearest
to Riemannian metrics, among all the geometric structures.
This is why, among
geometric dynamics, one might claim that the Lorentz dynamics is
the simplest one after the Riemannian dynamics. It is worth saying that
 this latter
dynamics is very poor. More precisely, let $f$
be an isometry of a compact Riemannian manifold $M$. Then
the closure
of an  orbit of $f$ is diffeomorphic to a torus, on
which $f$ induces a linear translation. One
can then say that there is only one pattern of
Riemannian dynamics: a linear translation on a torus.
On the other hand, one can  (roughly) characterize
Riemannian isometries, by topological
properties, such as  equicontinuity.

We know many patterns of Lorentz
dynamics. For instance, on the unit tangent bundle
of a hyperbolic compact surface, co-exist, the circle, the horocycle
 and the geodesic flows. All of them
preserve
the Killing Lorentz metric, although they present
antagonistic dynamics. Now the question is: find
(fundamental) patterns with the help of which, one
can build up any Lorentz isometry. Also, what are
the topological
properties
characterizing diffeomorphisms that preserve  Lorentz metrics?

Let's firstly observe that
   the Liapunov theory is not efficient
for our purpose. Indeed,
this theory is only sensitive
 to
 hyperbolicity,  since it
deals with the asymptotic
exponential behavior of vectors
 under the derivative of a  dynamical system.
For instance, this theory does  not
 distinguish
 between
the circle and horocycle  flows !

Here we will introduce
 another  growth notion, which appears to
be efficient
in studying some dynamical systems
with geometric properties (at least permitting to detect
a difference between the horocycle and circle flows !).

 Essentially, our fundamental growth notion is
the following:
Consider a diffeomorphism
$f$ of a compact manifold $M$.
We say that a vector $v \in TM$ is $f$-{\bf approximately
stable},
if $v$ is a limit of a sequence $(v_n )$ (in
 $TM$)
such that the image sequence $(Df^nv_n)$
 is bounded (i.e. lies in a compact subset of $TM$). We
denote
the set of such vectors by $AS(f)
\subset TM$ (observe that $AS(x, f) =
AS(f) \cap T_xM$ is not a priori a  vector subspace
of $T_xM$).

The  weakness of this
notion with regard to the Liapunov theory is clear,
the set $AS(f)$ may behave very badly.
Its efficiency
appears  in
 Lorentz dynamics. For a circle flow,
the situation
is trivial:  $AS = TM$. For the
horocycle
flow,  $AS$
 is the codimension 1 weakly stable or unstable
bundle of the geodesic flow, that contains the
given horocycle  flow (so $AS$ has dimension
2 and not 1). For the geodesic flow, $AS$ is exactly
its weak stable bundle.

The following theorem is one of our main results. It claims
that the codimension 1 property
of $AS$ is  a generic
situation in Lorentz dynamics. It seems also, that this
same property, gives a (rough) topological
characterization
of Lorentz isometries (however, we won't develop this
point of view here).

\begin{theorem} \label{theo.isometry} Let $M$ be a compact
Lorentz manifold (i.e. $M$ is endowed with
a  non degenerate symmetric 2-tensor
of signature $-+\ldots+$).
Let  $f $  be an isometry of $M$
 generating a
  non-equicontinuous subgroup
of isometries of $M$ (that is $\{ f^n / n \in {\bf Z} \}$
 is not pre-compact in  Isom$(M)$). Then the following
statements hold:

\noindent
{\bf Existence of approximately stable foliations.}
  The approximately stable set $AS(f)$
is a codimension one bundle tangent to
 a codimension 1 Lipschitz foliation
 ${\cal A S} (f)$, called the approximately
stable foliation of $f$. The leaves
of this foliation
are  geodesic and lightlike (that is the restriction
of the Lorentz metric to them is degenerate).
 The isometry $f$
preserves both its approximately stable foliation, and
 that  of $f^{-1}$, which  will be called
the approximately unstable foliation of $f$.

\noindent
{\bf Dynamics of $f$.} In the subset of points
where $AS(f)$ and $AS(f^{-1})$ are transverse, their
 (one dimensional) normal spaces $AS(f)^\perp$
and $AS(f^{-1})^\perp$ are respectively the
negative and positive Liapunov spaces of $f$.

The dynamics of $Df$ on the projective tangent bundle
${\bf P}(TM)$
 has the following ``north-south'' description:
Let
$U$ and $V$ be respectively neighborhoods of
the associated projective bundles ${\bf P} (AS(f))$
 and  ${\bf P} (AS(f^{-1}))$, then
$Df^n ({\bf P}(TM) - U)
\subset V$, for $n$  sufficiently large.

\noindent
{\bf Ergodic properties.}  Let $O$ be an open $f$-invariant
subset and  $\sigma$ a continuous  $f$-invariant
function defined in $O$. Then $\sigma$ is
constant along the leaves of the restriction
to $O$ of the two 1 dimensional
foliations tangent to $AS(f)^\perp$ and $AS(f^{-1})^\perp$.

 Finally,  the topological entropy of
$f$ vanishes exactly when the approximately stable and unstable
foliations of $f$
are identical.

\end{theorem}

It is a nuisance that we were not able to
exclude,   a priori,
 the approximately stable
and unstable foliations of an isometry,
may be different but coincide somewhere!


\paragraph{Generalized dynamical systems.}
The following terminology
will be useful.
 Consider a sequence $(f_n)$ of
diffeomorphisms of a compact manifold
$M$.
To emphasize  the fact  that
 we are  treating such a sequence
  in
a  dynamical viewpoint,
 we shall  call it
a {\bf {generalized dynamical
system}}.
For instance we  say that a vector $v \in TM$ is
 $(f_n)$-{ \bf {approximately
stable}},
if $v$ is a limit of a sequence $(v_n )$ in
 $TM$
such that the image sequence $(Df_nv_n)$
 is bounded.  We
denote
the set of such vectors by $AS((f_n))
\subset TM$.

So,  classical dynamical systems correspond
to the case: $f_n = f^n$,  where
$f$ is a diffeomorphism. But, even the less classical
case $f_n = f^{k_n}$, for a sequence
$( k_n )$ of integers (for example the return
times to some subset) is very interesting. Our philosophy
here is that, sometimes, one doesn't need the rich structure
of classical dynamical systems, especially the  various associated
cocycles.  Conversely, sometimes (in fact in many
situations in dynamical systems),  one can't avoid use (implicitly)
of generalized
dynamical systems. For instance, in our situation,
 it might happen that the theorem above is
dramatically empty, because every isometry
$f $ of $M$  generates a pre-compact subgroup, although
the group Isom$(M)$ itself is non-compact (for
example, Isom$(M)$ might be an infinite discrete torsion group).
In fact the proof of the theorem above, will
 come via the following generalized version.

\begin{theorem} \label{theo.generalized} Let $M$ be a compact
Lorentz manifold and
  $(f_n)$   a
 non-equicontinuous sequence
of isometries of $M$ (that is $(f_n)$
 is not contained in a compact subset of Isom$(M)$).

Then there is a  subsequence $(\phi_n)$
such that the approximately stable set $AS((\phi_n))$
is a codimension one bundle tangent to
 a codimension 1 Lipschitz foliation
 ${\cal A S} ((\phi_n))$, called the approximately
stable foliation of $(\phi_n)$. The leaves
of this foliation
are  geodesic and lightlike (that is the restriction
of the Lorentz metric to the leaves is degenerate).

After passing to a subsequence, we can assume that
the same is true for $( \phi^{-1}_n )$. In this case,
if $v \in TM-AS((\phi_n))$, $D\phi_nv $ tends to
$\infty$, and converges projectively (i.e. after normalization)
to  $AS{ (\phi_n^{-1})}$. The convergence is uniform
in compact subsets of $TM-AS((\phi_n))$. That is, if
$U$ and $V$ are  neighborhoods of
the associated projective bundles ${\bf P} (AS(( \phi_n)))$
 and  ${\bf P} (AS(( \phi_n^{-1})))$,  respectively,  then
$D\phi_n ({\bf P}(TM) - U)
\subset V$, for $n$  sufficiently large.

In addition,   we have the following
ergodic property.
Let $\sigma: O \to {\bf R}$ be a continuous
function defined in an  open
subset $O$. Suppose that
$\sigma$
 and $O$ are invariant
  by each $f_n$.
Then $\sigma$ is
constant along the leaves of the restriction
to $O$ of the two 1 dimensional
foliations tangent to $AS((\phi_n))^\perp$ and $AS((\phi_n^{-1}))^\perp$.

\end{theorem}

More details and a non compact variant of this theorem are
in \S \ref{lorentz.dynamics}.

\section{Results on isometry groups}

\subsection{Compactification}

The space of codimension 1 geodesic foliations
of a compact affine manifold  $M$ (i.e. $M$ is endowed with
a torsion free connection) is
naturally compact (this is because such a foliation
is uniformly Lipschitz, see \S \ref{lorentz.dynamics}). The same is true
for the space of codimension 1 geodesic
lightlike
foliations if $M$ is Lorentzian. We denote  this space
by  ${\cal F G}$. We construct a topology
on  Isom$M \cup {\cal F G}$, essentially
by the following rule: a sequence $( \phi_n )$
 in Isom$M$ converges to $F \in {\cal  F G}$
if $F$ is the approximately stable foliation of $( \phi^{-1}_n)$
 (that is if  $AS (( \phi^{-1}_n ))  = TF$).
One then shows that  Isom$M \cup {\cal F G}$
is a (metrizable) compact space, and that the action
of Isom$M$ on it is continuous.  Therefore,
by taking its  closure we get a natural compactification
of  Isom$M$. The same construction yields a compactification
and a {\bf boundary} $\partial_\infty \Gamma$ for any closed subgroup
of Isom$M$.

As example, in the case of the unit tangent bundle $PSL(2, {\bf R}) / \Gamma$
of a hyperbolic surface
${\bf H}^2/ \Gamma$, the isometry group is
$PSL(2, {\bf R})$.  Its boundary is
the circle, endowed with the usual projective action, see \S \ref{examples} for
further examples.

\subsection{Elementary groups}

The  ``north-south'' property described in
the two preceding theorems
tends
to suggest that the action of Isom$M$
on its boundary is of {\bf convergence} type, as in the case of
Fuschian groups \cite{Tuk}. This analogy between
Lorentz isometry groups and Fuschian groups
 may arise  in dynamical as well as in geometrical
 contexts.
In  \S \ref{bundle.compactification} we will see  how a Lorentz
geometry gives rise
to a fiberwise hyperbolic geometry, i.e. a family of
 hyperbolic spaces...
 This suggests that one could  translate
and update
notions on Fuschian groups to Lorentz isometry
 groups. Let's formulate
the simplest one.

\begin{definition} \label{def.polarized} A compact Lorentz manifold $M$
is called {\bf polarized} if the boundary of its isometry group
is one point.
\end{definition}

As example, an important geometrical class of polarized Lorentz manifolds
(that we won't consider here)  consists  of manifolds for
which the holonomy group
is reducible, but non-decomposable (this can't happen
 in Riemannian geometry).
 In this case the polarizing foliation is parallel, i.e. its tangent
bundle is invariant by parallel transport. (in this case the term
``polarization''
 may be  justified by physical, symplectical or
 optical considerations).

The Fuschian notion which  corresponds  to polarization is that of elementary
 parabolic groups. Now, elementary hyperbolic Fuschian groups
are those with limit sets of cardinality 2. In the Lorentz case,
it is not obvious  (at least at this stage) how  to check that
parabolic and hyperbolic behaviors don't co-exist. So, we were  able
to define  elementary groups only  via the following somewhat technical notion.

\begin{definition} \label{def.bifoliation} A codimension 1 lightlike
 geodesic {\bf bi-foliation} is
a map $x \in M \to L(x)=$
 a pair of
hyperplanes of $T_xM$, such that
there exist two codimension 1 lightlike geodesic foliations
${\cal F}^1$ and ${\cal F}^2$, with
$L(x) = \{ T_x{\cal F}^1_x, T_x{\cal F}^2_x \}$. We
will
denote such a bi-foliation by ${\cal F}^1 \cup
{\cal F}^2$.
\end{definition}

For
instance,  a pair  of two (perhaps identical) foliations
 gives a  bi-foliation.
 Notice,  however,  that in our definition, the
 foliations ${\cal F}^1$ and ${\cal F}^2$ are not part of the data.
Any foliation ${\cal F}$ such that $T_x{\cal F} \in L(x)$, for any
$x \in M$,  is called
tangent to
the given bi-foliation.
The nuisance here is that the coincidence
locus $C$ where the bi-foliation is one-valued, i.e.  where
 the two
defining foliations are tangent,  may be non-empty, or even worse,
 $M-C$ may have infinitely many connected components
(because one hasn't assumed things are analytic). In this case
the space of tangent
foliations (to the given bi-foliation) is an infinite compact space!
 Here, by analogy, not with Fuschian groups, but with higher rank
groups, one may call this space,
the apartment generated
by the bi-foliation, and denote it by  Apa$({\cal F}^1, {\cal F}^2)$
(in despite of the reference in notation to the foliations ${\cal F}^1$
 and ${\cal F}^2$, the apartment is clearly related to the
bi-foliation, only).

\begin{definition} \label{def.elementary} We say that a
 closed non compact subgroup
of Isom$M$ is {\bf elementary} if it preserves a codimension
1 lightlike geodesic bi-foliation. We say that $M$
is {\bf bi-polarized} if its isometry group is elementary.
\end{definition}

The following is a characterization of elementary groups.

\begin{theorem} \label{theo.elementary} Let $\Gamma$ be a closed
non-compact subgroup
of Isom$M$. Then $\Gamma$ is elementary if and only if
its  boundary has cardinality 1 or 2.

In the case $\partial_\infty \Gamma = 1$ point, there is no other
1-codimension geodesic lightlike foliation, preserved
by $\Gamma$. All the elements of $\Gamma$
have vanishing  topological entropy.

In the case $\partial_\infty \Gamma = 2$ points, up
to a finite index subgroup, $\Gamma$ is a direct
product of a compact group by
${\bf Z} $ or ${\bf R}$. Any element of the
${\bf Z}$ or the
${\bf R}$ factor has  positive entropy.

\end{theorem}

\begin{remark} {\em
In the case  where the boundary
of $\Gamma$ is
 is one point, we obtain nothing but the
Alexendroff compactification.   In fact  the group $\Gamma$
may be very large in this case.  For instance,  there are many
compact
homogeneous polarized Lorentz manifolds.
 Their isometry groups contain
solvable groups
obtained as semi-direct products of the circle
$S^1$ with Heisenberg groups (see \cite{A-S.1}
 and \cite{Zeg.espacetemps.homogenes}).
}
\end{remark}

\subsection{Amenable  groups. Dynamical structure of
 one parameter groups}

We also have the following
  convergence
 aspect
of the action of Isom$M$ on its boundary.

\begin{theorem} \label{theo.amenable} A closed non-compact
 amenable subgroup of
Isom$M$
is elementary.
\end{theorem}

Theorem \ref{theo.isometry} concerning groups isomorphic to
${\bf Z}$ will be a consequence of the two
above
theorems.
 Similarly, the  following result describes
the nice dynamical structure of non-equicontinuous
 isometric flows.

\begin{theorem} \label{theo.flow} Let $X$ be a Killing
field on a compact Lorentz manifold $M$,
generating a non equicontinuous  flow $f^t$.

Then $X$ is everywhere isotropic or spacelike, that
is $<X(x), X(x)> \geq 0$. Furthermore, $X$ is
contained in (and hence $f_t$ preserves)
two codimension 1
geodesic lightlike
foliations, called
 the approximately stable and approximately
unstable foliations of $f_t$.

These two
foliations  coincide  on the (invariant)
set of points where $X$ is isotropic (but the coincidence may
occur in a larger set).
 (In
particular,  if
$X$ is everywhere isotropic, then
 $X^\perp$ is integrable. This
is false for equicontinuous isotropic Killing fields).

\end{theorem}

\begin{remark} {\em
 It may  happen that
the approximately
stable and unstable foliations are identical, as
in the case of the horocycle flow. They might also,
a priori,  coincide in some proper subset of $M$. }
  \end{remark}

\subsection{Ergodic properties of
 non bi-polarized manifolds} \label{ergodic.tas}

Let $M$ be a compact non bi-polarized manifold,
 that is by definition, it has a non compact
and non elementary isometry group. Let
$S$ be the set of isotropic
vector fields $X$, which are tangent
to
foliations ${\cal F}$ belonging to
the boundary of Isom$(M)$ (note that because
$X$ is isotropic, then $X $ tangent to ${\cal F}$ is equivalent
to $X$ normal to ${\cal F}$). We denote
$S_x = \{X(x) / X \in S\}$ the evaluation
of $S$ at $x$, and
$TAS_x \subset T_xM$ the linear space that it generates.
 The  ``bundle'' $TAS$,  obtained in this
way will
be called the {\bf total approximately stable}
 bundle.
In a straightforward way, one sees that
somewhere the dimension of $TAS$ is $\geq 3$, because
 Isom$(M)$
is not elementary.
 In
the open
subset where  its dimension is
locally constant,
$TAS$ is a Lipschitz bundle. The same is true
for
its orthogonal $TAS^\perp$.
Observe that this  last bundle
is integrable. Indeed the geodesic
foliation obtained as the intersection
of all the foliations belonging to
$\partial_\infty$Isom$(M)$ is  a
natural candidate to be tangent to
$TAS^\perp$.  The starting point of Part II of this work, will
be the claim that  $TAS$ is also integrable,
in fact with umbilical leaves, which yields
 a local warped product structure for $M$. Here
we won't deal with
these geometric details but instead
with the following ergodic property
behind their proofs.

\begin{theorem} \label{theo.tas}  Let $E $ be a measurable linear subbundle
of $TM$,  containing $TAS$ and invariant by Isom$(M)$.
 Let $t: E \times E \to TM$
be a  measurable invariant bilinear bundle map.
Let ${\cal U}$ be the open set
$\{ x \in M /  \hbox{dim}TAS_x \geq 3 \}$. Then there
is a vector field $n$  defined almost everywhere in
 ${\cal U}$,  such that
$t(X, Y) = <X, Y>n$,  when $X \in TAS$
 and $Y \in E$.
 In addition,  we have the following ergodic property.
Let $\sigma: O \to {\bf R}$ be a
$C^1$ function defined in an  open set $O$,
and invariant by
 Isom$(M)$. Then $\sigma $ is constant
along $TAS$, that is $TAS$ is contained
in $Ker(d \sigma)$.

\end{theorem}

\subsection{Approximate isometries}

Questions on isometry groups of (fixed)  Lorentz
metrics on a manifold $M$, may be translated
 in a more stable way to
questions on the dynamics of the action of the
group of diffeomorphisms Diff$(M)$ on the space
of Lorentz metrics Lor$(M)$. So, the compactness of
Isom$(M, g)$ for any Lorentz metric $g$
on $M$, means that the stabilizers of the Diff$(M)$-action are compact.
A well known result of D'Ambra \cite{D'A} states this is true,
for compact simply connected manifolds, but, only in the analytic category.
 In fact
one  could  ask for the more stable:

\begin{question} Let $M$ be a compact simply connected manifold. Is
the action
of Diff$^2(M)$ on Lor$^2(M)$ proper ? (the exponent 2 stands for
 the $C^2$  differentiability  and the $C^2$ topology).

\end{question}

Also here the language of generalized dynamical systems is pertinent.

\begin{definition} Let $g$ be an element of   Lor$^2(M)$ and
$(f_i )  \subset $ Diff$^2(M)$ a generalized dynamical system.
 We say that $(f_i)$ is approximately isometric for
$(M, g)$  if $f_i^*g \to g$ in Lor$^2(M)$.
\end{definition}

Such dynamical systems occur exactly when the
action of
Diff$^2(M)$ on Lor$^2(M)$ is not proper.
In this article, we will not treat approximately isometric
dynamical systems
 (to avoid interference with
the approximate stability notion, which is  more
 central for us). However,
most of the results,
in particular Theorem \ref{theo.generalized},
generalizes, with the same proofs, to this context. As a corollary we have:

\begin{theorem} Let $M$ be a compact manifold such that
Diff$^2(M)$ acts non properly on Lor$^2(M)$. Then,
$M$ admits a codimension 1 lightlike geodesic foliation,
in the sense of some metric $g \in $ Lor$^2(M)$.
\end{theorem}

From a result of \cite{Zeg.geodesic.foliations}, a 3-manifold admitting
such a foliation is covered by ${\bf R}^3$. So, we have
the following partial  positive answer.

\begin{theorem} Let $M$ be a compact 3-manifold not covered by
${\bf R}^3$. Then Diff$^2(M)$ actsproperly
on Lor$^2(M)$.

\end{theorem}

\subsection{Organization of the article}

In the sequel, in order to simplify notations,
 we shall restrict our presentation
to
 compact manifolds. At certain points,
 we shall
remark on the differences in definitions and
 statements, for
 the non compact case.
 Our investigation of approximate stability
starts with the linear case
and is contained in  \S \ref{calculations}.  The general
affine
(i.e. connection preserving)
 appears in  \S \S \ref{affine.dynamics}, \ref{propagation}
 and \ref{theo.affine.existence}, and  is followed by
the Lorentz case in  \S \ref{lorentz.dynamics}.  Next , we introduce
a bundle compactification (\S \ref{bundle.compactification}), a
first step towards
the (foliation)
compactification of isometry groups of Lorentz manifolds
(\S \ref{foliation.compactification}). Theorem
\ref{theo.generalized} is proved in
\S \ref{lorentz.dynamics} and \S \ref{bundle.compactification}.
 In \S \S \ref{elementary} and \ref{partial.hyperbolicity}
we study elementary
 groups. This  part
 is somewhat technical due to the fact of   the possible co-existence
of parabolic and hyperbolic behavior for elementary groups.
 The last sections are devoted
to proofs of theorems
  and exposition of illustrative examples.

\paragraph{Acknowledgments.} I would like to thank M. Messaoudene
and G. McShane  for their helps.

\section{Approximate stability}

Here we give precise definitions
of   approximately  stable objects associated to
a {\bf {generalized dynamical system}}
 $(M, (f_n))$ (that is $(f_n)$ is just
a sequence of diffeomorphisms on $M$).

\begin{definition} \label{def.approximate}
 Let $(f_n)$ be a generalized dynamical system.

1. Stability.
 A sequence of vectors $(v_n)$
 in $TM$
is {\bf  stable} if
 $(f_n(v_n))$ is bounded.
  Similarly,
 a sequence $P_n \subset T_{x_n}M$ of tangent linear subspaces
 (of some fixed dimension)
is stable,
 if $D_{x_n}f_n \vert P_n$ is bounded.

2. Approximate stability.
 A vector $v \in T_xM$
is  {\bf {approximately  stable}} if $v$ is a limit of a stable
 sequence of vectors, that is there exists a convergent sequence
 $v_n  \to v$, such that
$Df_nv_n$ is  bounded. We denote by  $AS((f_n))$
the set of approximately stable vectors of $(f_n)$
and by $AS(x,  (f_n)) \subset T_xM$ its intersection
with $T_xM$.

We say that $v$ is {\bf { approximately
 strongly
stable}} if in addition  $Df_nv_n \to 0$ (in $TM$).

3. Punctual approximate stability.
 A vector  $v$ is  {\bf { punctually approximately
stable}}, if  it is a limit of stable sequence in the same
tangent space
$T_xM$: $v_n \in T_xM$, $v_n \to v$ and $D_xf_nv_n$ bounded.
As above we introduce  the notations
$PAS((f_n))$ and $PAS(x, (f_n))$.
We
 define analogously
a notion of {\bf { punctually approximately strongly
stable.}}

Similarly,  we define the same notions for tangent
subspaces in $TM$.

\end{definition}

\begin{remark} [Punctual and non punctual] {\em
Observe that $PAS(x, (f_n))$ is a linear space but  not
 $AS(x, (f_n))$ is. It is also important
to observe that if $P \subset  T_xM$ is such that
each $v \in P$ is approximately stable (i.e.
 $P \subset AS(x, (f_n))$), then
one cannot infer
that $P $ is an approximately
stable subspace. The vectors of
$P$ may be approximated by stable sequences
of vectors with different
base-points, and there is no way to sum (or
 even approximately sum)
them (since we have no control on the distance
between base-points). In the punctual case, it is clear
that $P \subset PAS(x, (f_n))$ is equivalent to
$P$ being  an approximately stable subspace. }
\end{remark}

\begin{remark}  [the non-compact case] {\em If $M$ is not compact, then
 in the definition of
stable sequence of vectors, we require that the sequence
$(f_n(x_n))$ lies in a compact subset of $M$.}
\end{remark}

\subsection{Sizes and modulus of stability} \label{def.modulus}

Although, when considering approximate stability,
we are far from uniform estimates,
we will sometimes  need  to
specify the
 sizes of objects. Since
 we are working on a compact manifold, we choose
a fixed continuous norm (on $TM$) with respect to
which, we
compare sizes.
In this case, we define
 the {\bf size} of a subset of $M$  to simply be  its
 diameter. Now,  let
$ P_n \subset T_{x_n}M $ be a stable sequence
 with respect to
 a generalized dynamical system $(M, (f_n)))$,  converging
to
 an approximately
stable subspace $P\subset T_xM$. By the {\bf modulus of  stability}
 of $P$ (or of $(P_n)$), we mean the inverse of the
 supremum of the norms
of $D_{x_n}f_n \vert P_n$.

However, to avoid explicit use of a particular norm, which would make
 the statements rather technical, and which has no natural
relation with data, we prefer to  thing of  sizes in a ``rough''
sense, that is up to multiplicative constants, given by the geometric
context.

For example,
given a stable
sequence $(P_n)$ (or even a family of such sequences
with uniform modulus of stability), one
can associate
 a continuous  family $(V_x)_{x \in M}$ where
$V_x$ is a neighborhood of
0 in $T_xM$, with size proportional
to that of the stable sequence. So,
if $P_n $ is a subspace of $T_{x_n}M$, then
 $Df_n(P_n \cap V_{x_n})$ is
contained in a neighborhood of $0$
in $T_{f_nx_n}M$, with uniformly  bounded size. For instance,  in  the
affine context  of \S \S  \ref{affine.dynamics}, \ref{propagation}
 and \ref{lorentz.dynamics},
we insist
 that the image by the exponential map
of the later neighborhood, is contained
in a convex neighborhood of
$f_nx_n$ in $M$.

\section{Examples. First calculations} \label{calculations}

\paragraph{The fundamental example.}
 We consider here (classical)
 linear dynamical systems
 on ${\bf R}^d$
 of the form $f_n = A^n$,
where $A \in GL(d, {\bf R})$.
 The fundamental example is
the following.
Let  $d=3$,
 endow  ${\bf R}^3$
with its canonical basis
$\{e_1, e_2, e_3\}$,
 and take:
$f_n = A^n=
\pmatrix{1&n&n^2/2 \cr 0&1&-n \cr 0&0&1}
=  \exp n \pmatrix{0&1&0 \cr 0&0&1 \cr 0&0&0}$.

Here approximate
stability and pointed approximate
stability coincide. At
$x= 0$, we have $AS(0) = PAS(0)= {\bf R}e_1 \bigoplus
 {\bf R}e_2$.
For instance to see that $e_2$
is
 approximately
stable, observe that the vector sequence $v_n =
(0, 1, -2/n)$ is stable. Indeed, $A^n(v_n) =
(0, 3,-2/n)$ is bounded.

Observe   that $e_1$ which is obviously stable, is
 in fact  strongly approximately stable.
To check  this consider $v_n = (1, 1/n^2, -2/n^2)$,  and note that
$A^n(v_n)= (0, 1/n^2+2/n, -2/n^2)$.

 Note,  however,   that in the case of ${\bf R}^2$
and $A_n = \pmatrix{1&n \cr 0&1 }
 = \exp n\pmatrix{0&1\cr 0&0}$, the approximately
stable space is reduced to ${\bf R}e_1$.

In general, consider  on  ${\bf R}^d$,
  $A_n = \exp n B$, where
$B$ is a Jordan block,
$B= \pmatrix{0&1&0&0.... \cr 0&0&1&0... \cr ...&0&0&1..
  \cr...&...&0 &0}$. Then
  $AS(0) = PAS(0)= {\bf R}e_1 \bigoplus
 {\bf R}e_2$, if $d \geq 3$.

Finally,  for $f_n = A^n$, with $A$
semi-simple, approximate stability
coincides with stability.

\subsection{Generalized linear
 systems  of ${\bf R}^d$ }

 In this section, we consider a
general sequence
 $(A_n) \in GL(d, {\bf R})$
of linear transformations of ${\bf R}^d$.

Let $A^+ \subset GL(d, {\bf R})$ denote
the semi-group
of diagonal matrices
diag$\{ \lambda_1, \ldots, \lambda_d\}$
where $0 <\lambda_1 \leq   \ldots \leq
\lambda_d$.
Then we can write $A_n = L_n D_n R_n$
where $D_n \in A^+$ and $L_n$ and $R_n$ belong
to $SO(d)$. Therefore,
a sequence of vectors  $(v_n)$
is stable in the sense of $(A_n)$
if and only if $(R^{-1}v_n)$
is stable in the sense of $(D_n)$.

Now,  it is straightforward
to calculate
the approximately stable
space for a
diagonal sequence. Indeed, as in the classical
case above,
 approximate
stability
 coincides with
stability. We point out  the following special case.

\begin{fact}
Let $A_n = L_n D_nR_n$ be such that $D_n$
has only one eigenvalue (with multiplicity)
$>1$, that is there exists $i$ such that
$\lambda_i \leq 1 < \lambda_{i+1}= \ldots =
\lambda_d$. Assume that $A_n$
is divergent  (i.e. has no convergent subsequence),
then a subspace $P $
is approximately
stable
 ifand only
there is a subsequence $Q_n$ of subspaces of ${\bf R}^i \times
 \{0\}$, such that
$P = \lim_{n \to \infty} R_n^{-1} Q_n$. In particular,
if $\lim_{n \to \infty}R_n^{-1} ({\bf R}^i \times \{0\})$ exists, then
it equals $AS(0, (A_n))$. Furthermore the moduli of
stability are always uniform (they  do not depend on
$(A_n)$).  In any case, we can always extract
a subsequence  of $A_n$, which has  an approximately stable space
of dimension $i$, and with uniform
modulus of stability.

\end{fact}

\paragraph{The linear Lorentz case.}
Here we consider sequences $(A_n)$ in
 $SO(1, d-1)$ i.e.   the orthogonal group
of the Lorentz quadratic form on ${\bf R}^d$:
 $q= -x_1^2+x_2^2+\ldots+x_d^2$. We note firstly
 that the    above observation holds.

\begin{fact}
An element  $A \in SO(1, d-1)$
can be written as $A = LDR$, with $L, R \in SO(d)$
and $D = $ diag$\{\lambda, 1, \ldots, 1, \lambda^{-1} \}$,
with $\lambda \leq 1$.

\end{fact}

\begin{proof}
 The $KA^+K$decomposition for  semi-simple Lie
groups yields in our case, $K = SO(1, d-1) \cap SO(d)
= SO(d-1)$, and $A^+$ any one parameter group
of symmetric matrices that
belongs to $SO(1, d-1)$. So, after conjugation by
a rotation $r \in SO(d)$, $rA^+r^{-1}$
becomes diagonal. The  eigenspace associated to an
 eigenvalue $<1$ (resp. eigenvalues $>1$) is
 isotropic (for the  form $q$).
 Hence,  it is one dimensional (because  $q$
is Lorentz). This proves the Fact.
\end{proof}

\begin{corollary} \label{linear.lorentz} Let $(A_n)$ be
a divergent sequence in
 $SO(1, d-1)$. Then,  there is a subsequence $(B_n)$
 such that the approximately
stable space $AS(0, (B_n))$ is a lightlike hyperplane.
The strongly approximately stable
space of $(B_n)$ is the orthogonal
of $AS(0, (B_n))$. It is an isotropic one-dimensional
space.
 The approximately
stable space
of $(A_n)$ is the intersection of
the hyperplanes  obtained from all the subsequences
 $(B_n)$. In particular,
 if
all the  approximately stable hyperplanes
 involved coincide, then, this equals the approximately
stable space of $(A_n)$ itself.

Finally, all the moduli of stability are uniform.
\end{corollary}

\paragraph{Discompactness.}
 In what  follows we give
an equivalent definition of the approximately
stable space, reminiscent to the Carri\`ere's notion of
 discompactness \cite{Car}. The codimension 1 fact
in the Lorentz case translates to
that $SO(1, d-1)$ has  discompactness 1. This
fact is as  crucial
for our work, as it  was  for  Carri\`ere's.

Consider
 a sequence $(A_n)$ in $GL(d, {\bf R})$,
and let $U$ be the unit ball of $R^d$. Then $E_n =
U \cap A_nU$ is a $d$-dimensional ellipse.
A limit (in the sense of Hausdorff)
of $(E_n)$ is an ellipse of dimension $\leq d$.
Let
$U^\prime$
 be the intersection of all the limits
of all convergent subsequences of $(E_n)$.
It is an ellipse of certain dimension.
Observe that $AS(0, (A_n))$ is the linear space
generated  by $U^\prime$.

\paragraph{Graphs.}

The approach of \cite{D-G} consists of taking graphs.  Keeping
the notation above,
consider the Graphs $Gr(A_n) =
\{ (x, A_n(x)    / x \in {\bf R}^d \} \subset
{\bf R}^d \times {\bf R}^d$.
Then $AS(0, (A_n))$ is  the intersection of all the projections
of all the limits of subsequences of $(Gr(A_n))$. In the Lorentz
case, we endow ${\bf R}^d \times {\bf R}^d$ with
the product $<,> \bigoplus-<,>$. Thus,  an element $A$ of
$GL(d, {\bf R})$ belongs to $SO(1, d-1)$ iff $Gr(A)$ is
isotropic. Observe that if for
some sequence $(A_n)$, $(Gr(A_n)) $
 converges to (a $d-$plane) $E \subset {\bf R}^d \times {\bf R}^d$,
then the intersection of $E$ with each of the factors ${\bf R}^d
\times \{0\}$ and $\{0\} \times {\bf R}^d$ are isotropic and  hence of
dimensions $\leq 1$. This is the content of discompactness 1.

\paragraph{Rank one groups.}
 The analogue of the fact above
is valid for
all the simple groups of non-compact type and of rank one, but
now we allow diagonal matrices of the form:
 diag$\{\lambda,
\ldots,\lambda, 1,\ldots, 1,\lambda^{-1},\ldots, \lambda^{-1}\}$.
 The multiplicity of $\lambda$
 (or $\lambda^{-1}$) may then be 1, 2, 4 or 8 (thanks to
the classification of simple Lie groups).

\paragraph{Chaos.}
Consider the two Lorentz linear systems:
$B_n=
\pmatrix{1&b_n&b_n^2/2 \cr 0&1&-b_n \cr 0&0&1}$
 and $C_n =
\pmatrix{c_n&0&0 \cr 0&1&0 \cr 0&0&c_n^{-1} }$. They
are orthogonal for
the form: $q = x_1x_3 + x_2^2$.
 We suppose that both $(b_n)$
and $(c_n)$
go to $\infty$ when $n \to \infty$. So we have
 as above,$AS(0, (B_n)) = {\bf R}e_1
\bigoplus
{\bf R}e_2$, and
 $AS(0, (C_n)) = {\bf R}e_2 \bigoplus{\bf R}e_3$.
 Let now $A_n = C_n B_n$ and observe
 that as above,
we have $AS(0, (A_n)) = {\bf R}e_1
\bigoplus
{\bf R}e_2$ and that $e_1$
is strongly approximately
stable, that is,  there is a sequence
$u_n \to e_1$ such that   $A_n(u_n) \to 0$.
 Of course,  $A_n(e_1) = c_ne_1 \to \infty$ !

\subsection{The derivative cocycle}

In the (Liapunov) measurable theory
one uses measurable trivializations
of $TM$, with respect to
which the derivative
of a diffeomorphism
$f$, is written as
a map $C_f: M \to GL(d)$, where $d = $ dim$M$.
In consideration of approximate
stability, one needs some control of the ``continuity''  of
the trivialization. So, to treat
the punctual approximate stability,
it suffices to consider a bounded trivialization.
That  is a frame-field
with image in a compact set of the frame
bundle. To handle
 approximate stability at some point, we
 further assume that the frame-field
is continuous at this point. So we will
always suppose that the trivialization
satisfy the needed requisitions .

Let now $(M, (f_n))$ be a generalized
dynamical system on a compact manifold
$M$. We denote $C_n = C_{f_n}$. In
the classical
case, i.e. $f_n = f^n$, we obtain
a cocycle $C: M \times {\bf Z} \to GL(d)$.
 In the generalized case, we just
obtain a collection of linear
systems
 $(C_n(x))$, for $x$ running over $M$.
 So one can relate the punctually approximately
stable space of $(f_n)$
at a point $x$ with the analogous one
of $(C_n(x))$ at 0.
\begin{fact} We have:
$PAS(x, (f_n)) = PAS(0, (C_n(x))) = AS(0, (C_n(x)))$.
\end{fact}

\paragraph{The Lorentz case.}
From the above facts we deduce:

\begin{proposition} \label{cocycle.lorentz}

Let $(f_n)$ be a non-equicontinuous
sequence of isometries
of a compact Lorentz
manifold $M$. Let $M^\prime$ be a countable
subset of $M$. Then there is
a subsequence $(\phi_n)$ of $(f_n)$
such that, for $x \in M^\prime$,
$PAS(x, (\phi_n))$ is a lightlike
hyperplane, and
$PAS(x, (\phi_n))^\perp = SPAS(x, (\phi_n)) $ ($=$
 the strongly
punctually approximately stable space
of $(\phi_n)$ at $x$). All these
hyperplanes have a uniform
modulus of stability.

\end{proposition}

\begin{proof} Keeping
the notations
above and using a Lorentz trivialization, we have at each $x$,
a derivative sequence
$(C_n(x) )$.  We shall see in \ref{start.unimodular},
since $(f_n)$ is not equicontinuous, that
for any $x$, $(C_n(x))$ is
not equicontinuous. The proof
follows by using  diagonal
procedure and  Corollary \ref{linear.lorentz}.
\end{proof}

\section{Affine dynamics: uniformity} \label{affine.dynamics}

Henceforth, we will only consider affine generalized
dynamical systems, that is $M$
 is  endowed with a linear
torsion free connection $\nabla$ and
$(f_n)$  is a sequence of  connection preserving
transformations.

\paragraph{Equicontinuity. Divergent sequences.}

The first fundamental property of
 affine dynamics is
the following:

\begin{proposition} [uniformity] \label{start.fact} Let $(f_n)$ be
 an affine
generalized dynamical system on a
 compact manifold $M$.
 Suppose that there is a
sequence $(x_n)$
such that $ ( D_{x_n}f_n )$ is
 equicontinuous
(that is $ ( D_{x_n}f_n )$ and
$ ( D_{x_n}f_n^{-1} = D_{f_n(x_n)}f_n^{-1} )$
 are  bounded). Then $ ( f_n )$ is
 equicontinuous,
that is,   $(f_n)$ lies in a compact
subset of the affine
group  Affin$(M)$.
\end{proposition}

The proof follows easily from
the fact that
 Affin$(M)$ acts properly on
the frame bundle
of $M$ \cite{Kob}.

\begin{corollary} [the unimodular case] \label{start.unimodular}
 Let $(f_n)$ be
 an affine unimodular
generalized dynamicalsystem on a
 compact manifold $M$, that is
all the $f_n$ preserve
a  volume form. (This is for example
the case if there is a parallel volume
form, e.g.  the connection
derives from
a pseudo-Riemannian metric).
 Suppose that there is a
sequence $(x_n)$
such that $ ( D_{x_n}f_n )$ is
bounded.
 Then $ ( f_n )$ is
 equicontinuous.
\end{corollary}
\begin{proof} The unimodularity and
the boundedness
of
$ ( D_{x_n}f_n )$, imply that
 $ ( D_{x_n}f_n^{-1} )$
 is also  bounded, and hence
 $ ( D_{x_n}f_n )$
 is equicontinuous.
\end{proof}

In the sequel we will be only interested in
the opposite situation of equicontinuity.
 Specifically, we
 say that a sequence $(f_n)$ is {\bf  {divergent}} if
$\{ f_n / n \in {\bf N} \}$ is a  closed discrete
subset of the group of homeomorphisms of $M$.
 So $(f_n)$ is not divergent if it contains  a
convergent subsequence (in the group
of homeomorphisms of $M$).

\paragraph{The codimension 0 case.}

We use the notations
of the  Proposition above.
 By definition $D_{x_n}f_n$ is bounded if
and only if $( T_{x_n}M)$ is a $(f_n)$-stable sequence
of subspaces (of codimension 0).
 In particular if $x_n \to x$, then
$T_xM$ is
an approximately stable subspace. So
 the corollary above
translates
to the fact that an unimodular
affine generalized
dynamical system
is equicontinuous
whenever $T_xM$
is an approximately subspace, for
some $x \in M$. So, in
particular
 a divergent sequence $(f_n)$ can
only have
approximately stable subspaces of
codimension
$>0$.

\paragraph{Modulus of stability.}
One may then ask whether or not
 the unimodularity
 condition
 is necessary. The answer is
yes, as we
shall see below in the case of Hopf
manifolds. These examples
also explain  why the phenomena of
``propagation of stability''
 (expressed
in \ref{propagation.stability}),
  which generalizes the uniformity fact \ref{start.fact},
 is only local (more precisely,  proportional to
modulus of stability). In fact,
the  pathology of Hopf manifolds is due
 to their
non-completeness. In affine flat dynamics,
 non-unimodularity and non-completeness
are generally thought of as being equivalent phenomena.

\paragraph{An example: Hopf manifolds. }
Recall that an affine (flat) Hopf manifold
is the quotient of ${\bf R}^d-\{0\}$
by a linear contraction. The simplest
case is when this contraction is
given  by a multiplication
map: $x \to \alpha x $, $0< \alpha <1$.
 The quotient $H_\alpha$ is thus endowed
with an affine action
of $GL(d, {\bf R})$. This action does not
 preserve any (non trivial)
measure. One can see this
in the following
 way.  For the sake of simplicity, we
only consider the
  the case $d=2$,
so that the Hopf manifold is
topologically a torus.
 Let
 $A^t$ be  a non compact
one parameter group of $SL(2, {\bf R})$.
 Then,   its orbits determine
 a Reeb foliation, and hence has only finitely
many recurrent leaves.

Let $f$ be the diffeomorphism
of $H_\alpha$
corresponding to a matrix
$\pmatrix{ 1 & 0 \cr 0 & \lambda}$
 with $\lambda >1$. Then $f^n$
corresponds to any of the matrices
 $ \pmatrix {  \alpha^{-m} & 0 \cr
0 & \lambda^n \alpha^{-m} }$, for
any integer $m$. Let $x \in H_\alpha$,
be the projection of a point $(a, b) \in
 {\bf R}^2 -\{0\}$.
 For calculation, one chooses a
fundamental domain containing
 $(a, b)$ and for each
$n$, one chooses
 $m$ such that
$ \pmatrix {  \alpha^{-m} & 0 \cr
0 & \lambda^n \alpha^{-m} } \pmatrix {a \cr b}$
belongs to the same domain.
 If $b = 0$, $x$ is a fixed point of $f$
 and so $D_xf^n$
is identified
with $ \pmatrix { 1 & 0 \cr
0 & \lambda^n }$.  Therefore,
 the approximately
stable space is the $x$-axis. For
$a \neq 0$, in
order to return to
the fundamental domain,
we will  need an integer
$m$ increasing with $n$
 and such that $\lambda^n \alpha^{-m} b$
is proportional to $b$.  Therefore,
$D_x f^n$ equals
$ \pmatrix {  \alpha^{-m} & 0 \cr
0 & \lambda^n \alpha^{-m} }$ and is thus equivalent
to $ \pmatrix {  \alpha^{-m} & 0 \cr
0 & 1/b }$ modulo multiplication by $\alpha$.   In particular,
$(D_xf^n)$ is bounded (but not equicontinuous).
Despite the fact the modulus
of stability is not uniform.

\section{Approximate stability in affine dynamics: partial uniformity}
\label{propagation}

In the present section, we prove some integrability and
geodesibility properties for the
approximately   stable objects, due to
a propagation of stability phenomenon, which holds  in
affine dynamics, and generalizing
the uniformity
Proposition \ref{start.fact}.

For $x \in M$, the exponential map $\exp_x $ is
 defined in an open subset
of $ {\cal D}ef_x \subset T_xM$.
 We recall that a submanifold $V$ of $M$ is geodesic
if whenever a geodesic $c: [a, b] \to M$ is
somewhere tangent to $M$: $c^\prime(t_0) \in T_{c(t_0)}V$
 then $c(t) \in M$, when $t$ belongs to some
neighborhood of $t_0$.
Notice that in general, if $P$ is a linear subspace of $T_xM$, then
$\exp_xP\cap {\cal D}ef_x$
 is not geodesic except  when dim$P= 1$
or for specific nice manifolds (see Part II).

\paragraph{Notations.}In the sequel,
for any $y \in M$,  we choose $V_y \subset T_yM$
 a neighborhood of $0$,  with size proportional
to the modulus of stability  of a given  $(f_n)$-stable sequence
of linear subspaces  $(P_n)$.

\subsection{Propagation of stability}

\begin{proposition} [Propagation of stability]
\label{propagation.stability} Let $( P_n)$
 be a $(f_n)$-stable sequence of
linear spaces, with $P_n   \subset T_{x_n}M$.
  Let  ${\cal P}_n = \exp_{x_n}(P_n \cap V_{x_n})$
 and $y_n \in {\cal P}_n$.
 Consider the restrictions $h_n = f_n \vert {\cal P}_n$. Then,
 the derivatives
$D_{y_n} h_n$  are
 uniformly bounded by   the size of
$(P_n)$.

 In particular, let $v_n
\in P_n$ be a convergent sequence of vectors,
 $y_n = \exp_{x_n}v_n$ and
$P^{\prime}_n = T_{y_n} (\exp_{x_n}P_n \cap V_{y_n})$. Then
$(P^\prime_n)$ is a $(f_n)$-stable sequence, with size
controlled by means of
that of $P_n$.
 \end{proposition}

\begin{proof}  Firstly,  observe that the claim is obvious
if $f_n$ are linear transformations of an Euclidean space
${\bf R}^d$, and $x_n = 0$ (here, without size restriction).

The proof of the general case follows by linearization.
  Assume to begin with that  the sequence
$(x_n)$ is stationary:  $x_n = x_0$,  and furthermore,
$x_0$ is fixed by all the $f_n$, $f_n(x_0) = x_0$, so that the
problem becomes linear
after  conjugation  by  $\exp_{x_0}$. More
precisely,  let $g_n = \exp_{x_0}^{-1} f_n \exp_{x_0}$ be defined on
some neighborhood $U$ of $0$ in $T_{x_0}M$.
It follows from the stability of
$(P_n)$, we may choose $U$ so that $g_n(P_n \cap U)$ is
contained in some fixed small neighborhood $U^\prime$ of $0$.
Hence,  the derivatives of $f_n = \exp_{x_0} g_n \exp_{x_0}^{-1}$
along points of $\exp_{x_0}(P_n \cap U)$ are comparable to the
corresponding ones for $g_n$,
 because of
the fact that  all things stay in a compact set, where
the derivatives of $\exp_{x_0}$
and $\exp_{x_0}^{-1}$ are controlled.

When  $(x_n)$ is not stationary,
 we consider the family of derivatives $g_n = D_{x_n}f_n: T_{x_n}M
 \to T_{f_n(x_n)}M$.
 Inasmuch these spaces
are equipped with norms induced from a metric on $TM$, with
respect to which we define stability, these norms are defined up to
a bounded distortion.
Now,  since by definition $x_n$ and $f_n(x_n)$ stay in a compact set, we
can find identifications of bounded distortion
of all our linear tangent spaces with
a fixed Euclidean space. Therefore the stability notions
are preserved
and the proof goes as in the previous case.

\end{proof}

\begin{proposition} [compatibility with parallel transport]
\label{propagation.parallel}
Keeping  the notations of the proposition
above, let $(c_n)$ be
a sequence of curves, such that the image
of  $c_n : [0, 1] \to M$ is
contained in
$\exp_{x_n} P_n \cap V_{x_n}$ with $c_n(0) = x_n$.
 Consider
  $P^{\prime \prime}_n$ the
 parallel
transport of $P_n$ along $c_n$. Suppose that $(c_n)$
 is bounded in the $C^1$
topology.
 Then $P^{\prime \prime}_n$
is a stable sequence.
\end{proposition}

\begin{proof}
From the proposition
 above the image curves $d_n = f_n(c_n)$ are
bounded in the $C^1$ topology.
Let $\tau_n = T_{c_n(0)} M \to
T_{c_n(1)}M$ be
 the parallel transport
along $c_n$ and let $\tau_n^\prime$ be the
 analogous parallel
transport along
$d_n$. Then,  $\tau_n$ and $\tau_n^\prime$
are uniformly bounded
(since $c_n$ and $d_n$ are $C^1$ bounded). Now
 because $f_n$
are affine, they commute with parallel
transport, in
particular: $D_{c_n(1)}f_n =
(\tau_n^\prime)^ D_{c_n(0)}f_n
 \tau_n^{-1}$. Therefore,  $(\tau_n (P_n))$ is a
stable sequence.
\end{proof}

\subsection{Geodesibility properties}

\begin{proposition}
 \label{maximal}  Suppose that
the stable sequence $(P_n)$ converges to
 an approximately
stable
 linear subspace  $P \subset T_xM$, which  is maximal among such
 subspaces (that
is not strictly contained in another approximately stable
subspace). Then $\exp_xP \cap V_x$
is geodesic in $M$. The same result is true in
the punctual approximately stable case, that is,  for any $x \in
M$, $\exp_x(PAS(x, (f_n)) \cap V_x)$ is geodesic.
\end{proposition}

\begin{proof}
  Let
$c_n: [0, 1] \to \exp_{x_n} P_n \cap V_{x_n}$ be
a sequence
of curves converging in the $C^1$
topology to a curve $c : [0, 1] \to \exp_xP \cap V_x$, and
 such
 that $c_n(0) = x_n$. Let $y_n = c_n(1)$,
$P_n^\prime \subset T_{y_n}M $ be the tangent
 space  of $\exp_{x_n}P_n \cap V_{x_n}$ at $y_n$,
and $P_n^{\prime \prime} \subset T_{y_n}M$  the parallel transport
of $P_n $ along $c_n$. By the above propositions, both
$(P_n^\prime )$ and $(P_n^{\prime \prime} )$ are
stable sequences of linear subspaces.

Denote the analogous objects at $y = c(1)$
by $P^\prime$ and $P^{\prime \prime}$, which are obviously the limits
of $P_n^{\prime}$
and $P_n^{\prime \prime}$,  respectively.

To prove that
$\exp_xP \cap V_x$ is geodesic, it suffices to check
the equality: $P^\prime = P^{\prime \prime}$. Indeed,  if
this is true for arbitrary $c$, then the tangent space of
$\exp_xP\cap V_x$ is parallel (along itself), which is
equivalent to the geodesibility (for torsion free connections).

But,  if $P^{\prime} \neq P^{\prime \prime}$,  then
$P^\prime \bigoplus
P^{\prime \prime}$
is an approximately stable linear subspace. Indeed,
let $\{e^i \}$ be a basis of this later subspace, and
choose $\{e_n^i\}$ vectors of
$P_n^\prime \bigoplus
P_n^{\prime \prime}$, such that $e_n^i \to e^i$.
Denote by $E_n$ the vector
subspace generated by the $\{e_n^i\}$.
Then  $(E_n)$ is a
stable sequence converging to $ P^\prime \bigoplus
P^{\prime \prime}$.
  This contradicts the fact
that $P$  is maximal. Therefore $P^\prime = P^{\prime \prime}$.

\end{proof}

\subsection{The codimension 1 case: lamination properties}

\begin{fact} \label{hyperplan}
Assume that $(f_n)$ has no codimension 0 approximately
stable subspace (that is for no $x \in M$, $T_xM$ is
an approximately stable space).
 Let $P \subset T_xM$ be an  approximately stable
hyperplane. Then $AS(x, (f_n))= P$.
 In particular, in this case $P$ is  maximal and hence from
the Proposition above, $\exp_xP \cap V_x$
is a geodesic hypersurface.
\end{fact}

\begin{proof}
Assume the contrary, that is there exists an approximately stable
 vector $v \in T_xM$ which is transverse to $P$. Thus,  $P$
 and $v$ are respectively limits of stable sequences $P_n
\subset T_{x_n}M$, and $v_n \in T_{y_n}M$. By transversality
(of $P$ and $v$), $\exp_{x_n}P_n \cap V_{x_n}$ and $\exp_{y_n} {\bf R}v_n
\cap V_{y_n}$
intersect, in an uniform transverse meaner at  a point
$z_n$ (near $x$, for $n$ large). Moreover by  Proposition
\ref{propagation.stability}, $D_{z_n}f_n$ is uniformly bounded along the
tangent spaces at $z_n$,  of each of the submanifolds
$\exp_{x_n}P_n\cap V_{x_n}$ and $\exp_{y_n} {\bf R}v_n \cap V_{y_n}$. By
definite
transversality, $D_{z_n}f_n$ is bounded in $T_{z_n}M$,
that is $(T_{z_n}M)$ is
a $(f_n)$-stable sequence,  a fact that  is excluded by hypothesis.

\end{proof}

\begin{corollary} \label{hyperplan.geodesic}

With the same hypothesis,
 if $y \in exp_xP \cap V_x $, then $AS(y) = T_y(\exp_xP\cap V_x)$.
\end{corollary}

\begin{proof}
  It follows from   Proposition \ref{propagation.stability},
 that $T_y (\exp_xP\cap V_x) \subset AS(y)$, and
we infer
 from the fact above, the  equality: $AS(y) = T_y(\exp_xP)$, as desired .

\end{proof}

By the same argument we get:
\begin{corollary} \label{disjoint}
  If $P_x \subset T_xM$ and $P_y \subset T_xM$ are
approximately stable hyperplanes, then the
 two geodesic hypersurfaces $\exp_xP_x\cap V_x$ and
$\exp_yP_y \cap V_y$ are either disjoint or tangent, in which
case their intersection is open in each of them.
\end{corollary}

\paragraph{The method of Graphs.} As in the linear case
 (\S \ref{calculations}),  there is a graph approach leading to
some geometric proofs of the above
 properties of the approximately stable
submanifolds $\exp_xP \cap V_x$. To see this, note that
if we endow w
the product
$M \times M$ with the product connection. Then $(\hbox{Graph}(f_n))$
 is a sequence of geodesic submanifolds
 in $M \times M$. Any approximately stable submanifold
$\exp_x P \cap V_x$
is obtained as a projection of a limit
 of a sequence of connected components of
$U \cap \hbox{Graph}(f_n)$, where $U$ is an open subset
of $M \times M$.

The geodesic character of $\exp_xP \cap V_x$
is thus obvious. Nevertheless, neither the lamination
properties (\ref{hyperplan}), nor the
control of sizes by modulus of stability, seems
to be easy to  treat, via
this approach.

\section{The approximately stable foliation theorem}
\label{theo.affine.existence}
Here follows a fundamental existence and regularity result:

\begin{theorem} \label{theo.existence}

Suppose that there is a dense subset
$M^\prime \subset M$, in which
$ PAS(x, (f_n))$ is a hyperplane,
with uniform modulus of stability. Then
$AS((f_n))$ is a Lipschitz
 codimension 1 subbundle of $TM$, tangent
to a geodesic foliation ${\cal A S}((f_n))$,
called the approximately
stable foliation of $(f_n)$.
\end{theorem}

For the proof we need:

\paragraph{Digression: Codimension 1 geodesic laminations.}

Geodesic laminations (and in particular foliations)
of codimension 1, enjoy some remarkable regularity
properties.
These  properties are well
known for laminations of hyperbolic
 surfaces \cite{Thu}, but are
in fact valid in the general
context  of
codimension 1 geodesic
laminations, in the sense of a
connection (and also some
some codimension 1 foliations with
other geometric
origins).  For proofs and
 related
questions,  see
\cite{Zeg.geodesic.foliations} and \cite{Sol}. Behind all
of the regularity
results  for  geodesic
laminations, is the following fundamental
Lipschitz regularity fact.

\begin{lemma} \label{geodesic.foliation.lemma}

 Let $M$ be a compact
manifold
endowed
with a torsion free connection and an auxiliary
norm $\vert.\vert$ on $TM$.
   Let $M^\prime$ be a subset of $M$
and suppose given a real $r$ and
for $x \in M^\prime$
 a hyperplane $P_x \subset T_xM$  and let
 ${\cal L}_{x, r} = \exp_x P_x \cap B_x(r)$
 where  $B_x(r)$ is  the ball of $T_xM$
centered at $0$ and with
 radius $r$.  Also, suppose that ${\cal L}_{x, r}$
is geodesic and that if
 two plaques
 ${\cal L}_{x, r}$
 and  ${\cal L}_{y, r}$ intersect at some
point, then they are tangent at that point
(and hence by geodesibility,
 the intersection  ${\cal L}_{x, r} \cap{\cal L}_{y, r}$
 is open in both  ${\cal L}_{x, r}$ and  ${\cal L}_{y, r}$).
 Then,  along $M^\prime$, the map $x \to P_x$
is Lipschitz, with Lipschitz
constant depending
only on the geometry of the connection,
the auxiliary norm and $r$.
\end{lemma}

We infer from the above lemma the following corollary:

\begin{corollary} \label{geodesic.foliation.corollary}

(i)  A codimension 1 geodesic lamination on $M$ is
 Lipschitz,  with Lipschitz constant depending
only on the geometry of the connection
and the auxiliary norm.
It then follows that
the space    of codimension 1 geodesic foliations
 of $M$, endowed with the $C^0$ topology (or equivalently
the Lipschitz topology) on hyperplane fields, is compact.

(ii)
With the same hypothesis as
  in the lemma above, suppose that the
set $M^\prime$
is dense. Then,  the geodesic plaques ${\cal L}_{x, r}$
extend to a geodesic foliation of $M$.
\end{corollary}

\paragraph{Proof of Theorem \protect\ref{theo.existence}}

 Since the modulus
of stability is uniform,
the geodesic plaques
$\exp_x PAS(x, (f_n))\cap V_x$
 given by \ref{hyperplan},
satisfy the conditions of \ref{geodesic.foliation.lemma}. Therefore,
 we have a  geodesic foliation
${\cal F}$ of $M$
such that $T_x{\cal F} = PAS(x, (f_n))$
 for $x \in M^\prime$. Again by
uniformity of the modulus of stability, and the assumption that
$M^\prime$ is dense, for any $x \in M$,
 $T_x {\cal F}$ is an
 approximately
stable hyperplane at $x$. Therefore, it follows from Fact
   \ref{hyperplan} that  for any $x \in M$,
$AS(x, (f_n)) = T_x {\cal F}$.
 $\Box$

\section{Lorentz dynamics. Proof of the existence part in
Theorem \protect\ref{theo.generalized} }
\label{lorentz.dynamics}

Henceforth, we will  only deal with
 transformations preserving a Lorentz
structure. So $M$ is now a compact Lorentz
manifold and $f_n \in $ Isom$M$.
 In the present section,  we shall
prove the existence part of
Theorem \ref{theo.generalized} that is
the existence
of approximately
stable foliations for subsequences (cf.  Fact \ref{fact.generalized}). The
remaining part of
Theorem \ref{theo.generalized}, will
be proved in \S \ref{last.part}. To begin with, we recall some facts
about  lightlike geodesic
foliations.

\paragraph{Lightlike geodesic foliations. The
 compact space ${\cal F G}$.}

A (codimension 1) foliation
$\cal F$ is {\bf {lightlike}} if
the restriction of the metric
to $T{\cal F}$
is degenerate. This  metric is
thus positive non-definite
and its Kernel is  a 1-dimensional
sub-foliation of $\cal F$, denoted by  $Nor{\cal F}$
and called  the {\bf  { normal foliation }}
of ${\cal F}$. This 1-dimensional foliation
is isotropic (i.e. the metric vanishes along it)
and completely determines $\cal F$, since
$T{\cal F}$ is just the orthogonal
of $T(Nor{\cal F})$.

Here,  we consider  lightlike geodesic foliations
\cite{Zeg.geodesic.foliations}. One
of their basic properties is that their (1-dimensional) normal
foliations are also geodesic.

We denote by ${\cal F G}$ the space of all  lightlike geodesic
foliations on $M$. It is a closed subset in the space
of all the geodesic foliations, and so it is compact whenever
$M$ is.

\paragraph{Existence of approximately stable lightlike
geodesic foliations.}

The existence part in  Theorem \ref{theo.generalized}
is a consequence of:

\begin{fact} \label{fact.generalized}
Let $M$ be a compact
Lorentz manifold
 and $(f_n)$ a non-equicontinuous  sequence of isometries.
Then,  there is a  subsequence $\phi_n$ admitting an
 approximately
 stable lightlike (codimension 1) geodesic foliation
${\cal A S}((\phi_n))$,
 that is the approximately stable set $AS((\phi_n))$
is a codimension one bundle tangent   to
    a codimension 1 geodesic lightlike  foliation.

 In fact,
 $AS((\phi_n))$ coincide with  $PAS((\phi_n))$,
the punctually approximately stable bundle.
Furthermore,
 the 1-dimensional normal foliation of ${\cal A S}
(( \phi_n ))$ equals the (punctually or not)
approximately strongly
 stable  bundle of $(\phi_n)$.
\end{fact}

\begin{proof} Choose
a countable dense subset $M^\prime$
of $M$. Then, by \ref{cocycle.lorentz}, there
is a subsequence $( \phi_n)$
 such that $PAS(x, (\phi_n))$
is a lightlike
 hyperplane, with
uniform
modulus
of stability, when $x \in M^\prime$. We may
infer from Theorem \ref{theo.existence},
the existence on $M$  of a
codimension 1 lightlike
geodesic foliation
${\cal A S }((\phi_n))$ tangent
to $AS((\phi_n))$.

On $M^\prime$, we have the equality:
$AS(x, (\phi_n)) = PAS(x, (\phi_n))$.
 To prove the equality for an
arbitrary $x \in M$, we
apply \ref{cocycle.lorentz} for $( \phi_n )$
(instead of $(f_n)$) and
 for $M^\prime = \{x\}$. We
obtain
a subsequence $ (\psi_n )$
 for  which
$PAS(x, (\psi_n))$
 is a hyperplane and
hence
equals $AS(x, (\phi_n))$
(since
obviously $AS(x, (\phi_n)) \subset
AS(x, (\psi_n))$
and $AS(x, (\psi_n)) =
PAS(x, (\psi_n))$
from \ref{hyperplan}).
 In particular
$PAS(x, (\psi_n))$ does not depend
on the subsequence $(\psi_n )$
(provided it is a hyperplane).
By \ref{linear.lorentz}, this
implies that $PAS(x, (\phi_n))$
 is itself  a hyperplane, and thus
equals $AS(x, (\phi_n))$.

 The fact that  the normal
direction of ${\cal A S}((\phi_n )$
  is  strongly approximately stable
 follows from the analogous
statement in the linear case, \ref{linear.lorentz}. In the same
fashion, one can prove
the coincidence
 between  strong approximate stability
and  strong punctual  approximate stability.
 \end{proof}

\paragraph{The non compact  case.}

If $M$ is not  compact, it may
happen that for some $x$, $AS(x, ( \phi_n)) = \{0\}$
 for any subsequence $( \phi_n )$ of $(f_n)$. Indeed, from
 our
definition \ref{def.approximate}, this
 happens when $x$  satisfies the uniform
escaping
property: for any sequence $x_n \to x$ and
a subsequence $( \phi_n)$ of $(f_n) $,
$(\phi_n(x_n) )$ tends
 to $\infty$ (i.e. leaves all compact
subsets of $M$).
 \paragraph{``Recurrence'' for generalized dynamical systems.}
There is no natural way to define recurrence
notions for generalized dynamical systems,  so
that they satisfy
a kind of Poincar\'e recurrence Lemma. For our purpose, the following
notion seems interesting:

\begin{definition} \label{def.non.escaping}

Let   $K$  be
a compact
subset of $M$. Define its  non-escaping subset
$NE(K, (f_n))$ as the (compact) subset of
points  $ x \in K $ such that
there is $x_n \in K$,
$x_n \to x$ and $f_nx_n \in K $.
In other words $NE(K, (f_n))$ is
the Hausdorff  limit of the sequence of compact sets
$ K \cap f_n^{-1}K$.
In particular $\hbox{Vol}(NE(K, (f_n))
\geq \hbox{Vol}K \cap f_n^{-1}K) \geq
2\hbox{Vol}(K)-\hbox{Vol}(M)$ (because the $f_n$
 are volume preserving). (In particular, letting
$K$
having a large relative volume, we see that almost every
$x \in M$
is non-escaping for some compact $K$).

\end{definition}

We lose uniformity when $M$
is not compact, and so to  estimate
sizes, we will choose a norm
$\vert.\vert$ on $TM$.
 As sizes may depend upon the  choice of
such a norm,
the previous uniformity of modulus of
stability (\ref{cocycle.lorentz}) fails.  A straightforward localization
of the previous arguments yields:

\begin{proposition} Let $M$ be a (not necessarily compact)
 Lorentz
manifold, $(f_n)$ a sequence in  Isom$(M)$ and  $K$
a  compact
subset of $M$.  Then there is a subsequence $( \phi_n )$
 such that,
on $NE(K, (f_n))$, $AS(( \phi_n ))$ is
a Lipschitz codimension 1 bundle, and there is
$r=  r(K)$ (not depending  on
$(f_n)$), such that
the family of plaques
${\cal L}_{x,r} = \exp_x AS(\phi_n)) \cap B_x(r)$ determine
a  codimension 1 lightlike geodesic lamination of a
neighborhood of $NE(K, (f_n))$, tangent
to $AS(( \phi_n))$ in $NE(K,  (f_n))$.
\end{proposition}

This applies at least in the finite volume
case, because for instance if
 Vol$(K) > $ Vol$(M)-\epsilon$, then,
  Vol$(NE(K, (f_n))) >$ Vol$(M)-2\epsilon$.
 Therefore, approximately  stable  bundles
of some subsequences of $(f_n )$, give rise to
 laminations, along
big volume subsets of $M$. However, in order to obtain
 a foliation, on the whole of $M$, we
must  find volume exhausting sequence of compact sets
$(K_j)$
 with non collapsing radii $(r(K_j) )$.
This is generally impossible because of
non compatibility of the auxiliary metric $\vert. \vert$
with  the natural data.

Therefore,
 to avoid use of such a norm,  we introduce
for $x \in M$,
 ${\cal D}ef_x \subset T_xM$ to be the domain of definition
of the exponential map $\exp_x$.
It is open (by definition) and star shaped in $T_xM$, that is,
if $u \in
{\cal D}ef_x $, then $tu \in {\cal D}ef_x$, for
$t \in [0, 1]$).  Consider also the
regular domain of definition   ${\cal D}ef^*_x$
 defined as the (open) set of vectors
$u \in {\cal D}ef_x$, such that
$tu$ is a regular point of $\exp_x$
for any $t \in [0, 1]$. Thus, if
$E $ is  a linear subspace of $T_xM$, then
$\exp_x(E \cap {\cal D}ef^*_x)$
is  an immersed submanifold of  $M$.
 This defines open subsets ${\cal D}ef$ and
${\cal D}ef^*$ of $TM$, which are
invariant by  Isom$(M)$.

The propagation of stability \ref{propagation.stability} is valid  a
priori only
in a domain having a size proportional
to the modulus of stability. The obstruction lies in
fact essentially in
non completeness, as this was  shown in
the examples of Hopf manifolds.  Therefore,  one may
 hope for
a propagation of stability
in ${\cal D}ef$. Indeed,  this was shown in
\cite{Zeg.actions.affines} in the context of (strict) stability (instead
of approximate stability here), but only
 generically. The point is that
 we just need to ensure continuity
 (by using the Lusin theorem)
of the map $x \to {\cal D}ef^*_x$, and also
the continuity of
its variants obtained by intersection with
stable (here approximately stable) bundles. The same
proof yields the following fact:

\begin{proposition} Let $M$ be a (not necessarily compact)
 Lorentz
manifold, $(f_n)$ a sequence of  Isom$(M)$.
 Let $\mu$ be a Borel finite  measure (not
related to data, but serves in applying the Lusin
almost everywhere continuity theorem). Given
a positive real $\epsilon$, there is
a compact
$K$
in
 $M$, with $\mu(K) > \mu(M)-\epsilon$ and
there is a subsequence $( \phi_n )$
 such that,
the family of plaques
${\cal L}_x = \exp_x AS( (\phi_n )) \cap {\cal D}ef^*_x$,
 for $x \in NE(K, ( f_n ))$,  determines
a  codimension 1 lightlike geodesic lamination of a
neighborhood of $NE(K, (f_n))$, tangent
to $AS( ( \phi_n))$ in $NE(K,  (f_n))$ (that is
if $y \in NE(K, (f_n )) \cap {\cal L}_x$, then
$AS(y, ( \phi_n )) = T_y {\cal L}_x$).
\end{proposition}

In the finite volume case, we use  this proposition for
the Lorentz measure. As we remarked earlier,
we get  Vol$(NE(K, ( f_n ) > $
Vol$(M) - 2 \epsilon$. By letting $K$ get  larger
and larger, and applying a diagonal
process, we obtain a foliation:

\begin{theorem} [The finite volume case] The Theorem
 \ref{theo.generalized}, extends to the case of  $M$ of finite volume: a
 divergent  sequence
$(f_n )$ of Isom$(M)$ possesses a subsequence
$( \phi_n )$ admitting an approximately
 stable foliation ${\cal A S}(( \phi_n ))$.
Furthermore, almost everywhere:
 $AS(x, ( \phi_n )) = PAS(x, ( \phi_n ))$ (this equality
is everywhere true in the compact case).
\end{theorem}


\paragraph{A relative  version.} Here we give a relative version:

\begin{theorem} Let $M$ be a finite volume
Lorentz manifold,
$(f_n)$ a divergent sequence of isometries
that leaves invariant
a finite volume Lorentz submanifold $N$.
Let $(\phi_n)$
be a  subsequence of $(f_n)$ having
an approximately stable foliation ${\cal A S } ((\phi))$.
 Then,
 the restriction $ ( \phi_n\vert N )$ admits
an approximately stable foliation, which is
just the trace of ${\cal A S}(( \phi_n ))$
on $N$. In particular,   the intersection of a leaf
of ${\cal A S} (( \phi_n ))$ with
$N$ is geodesic in $N$.
\end{theorem}
This result   follows from the fact that: if $x \in N$, then
 $PAS(x, ( \phi_n \vert N )) = PAS(x, ( \phi_n )) \cap T_xN $

\paragraph{A foliated version.}
We also have the following foliated version:
\begin{theorem}  \label{theo.foliated.version}
Let $M$ be a finite volume
Lorentz manifold
 and $(f_n)$ a divergent sequence of isometries.
Let $U $ be an open subset of $M$
and
 $L \subset TM\vert U$  a  $C^0$ subbundle over $U$, invariant
by each $Df_n$ (in particular $U$ itself  is invariant by
 $f_n$). Let $M_0$ be the subset of $U$ where
$L$ is of lorentzian type (i.e. the restriction
of the metric has signature $-+\ldots+$).
Let $(\phi_n)$
be a  subsequence of $(f_n)$ having
an approximately stable foliation ${\cal A S } ((\phi))$. Then
 in $U-M_0$, $L \subset  AS((\phi_n))$
and in
$M_0$, $L^\perp \subset AS((\phi_n))$ (where
$L^\perp$ is the orthogonal of $L$).

Suppose that $L$ is integrable. Then,
 in $M_0$, the leaves of
${\cal A S} (( \phi_n )) \cap L$
 are geodesic inside
the leaves
of $L$.
\end{theorem}




\section{Bundle compactification} \label{bundle.compactification}

Let $M$ be a manifold.  One
may naturally construct a fiber metric
on its (principal) frame bundle $P_M \to M$,
 so that any
 $C^1$ diffeomorphism acts isometrically
on the fibers.  In fact,  instead of $P_M$,
it is
more convenient to consider the associated
bundle $S_M \to M$
whose typical fiber
the universal symmetric
Riemannian space $S_k =
SL(k)/SO(k)$, where $k = $ dim$M$. So, $S_M$ may be interpreted
as the bundle
of conformal structures on (the fibers of) $TM$.

 One may then naturally  compactify
the fibers to get a bundle
$\overline {S_M} \to M$ with a compact fiber
$\overline {S_k}$, the Hadamard compactification
of $S_k$ which is topologically  a closed ball.

It is also interesting to interpret
 recurrence properties of the
action of $f$, or more generally,
 of a generalized dynamical system
 $(f_n)$  on $\overline{S_M}$.  The philosophy is that
recurrence conditions may be related to
 classical
notions  as,
 Oseledec's decomposition, invariant metrics (of some
regularity)...In the sequel we will rather use this construction  to
interpret the approximate stability.

Notice that  these
constructions generalize to  fiber dynamical systems
on a linear bundle   $E \to M$.
 It is also possible to  consider the case of bundles with a $G$-structure
(i.e a reduction of the structural group
 of the principal bundle $P_M \to M$ to $G$) where
$G$ is a semi-simple subgroup
of $SL(k, {\bf R})$ of non compact type.
 In which case, one  can build the same constructions with the symmetric
space
associated to $G$.  In fact,  instead of developing
the general situation, we shall henceforth restrict  ourselves to
 the case where
$G $ is $ SO(1, k-1)$, and  $E$ is the tangent bundle
of $M$. This exactly means that
$M$ has a Lorentz structure. We mention at this stage that
some but not all of the
next results   are valid in the
general case.

\paragraph{Fiberwise hyperbolic geometry.}

Let $(M^{1+d}, <,>)$ be a Lorentz manifold.
 Note $T^rM = \{v \in TM
 / <v, v> = r \}$.
 After passing if necessary to a double cover of $M$,
we can choose a sheet of
 $T^{-1}M$  that we note  ${\bf H}M$. It is a bundle over
 $M$  with type fiber the hyperbolic space
 ${\bf H}^{d}$ (recall that the Lorentz metric
has signature $-+\ldots+$). It is compactified
by adding the
the  projectivization  ${\bf S}^\infty M$
 of the isotropic cone  $T^0M$.
 We denote it by
$\overline{{\bf H}}M$ the bundle over $M$
with type fiber a topological
closed ball of dimension $d$.

The group  $Isom(M)$
acts  on ${\bf H}M$
and on  $\overline{ {\bf H}}M$,
 by preserving  the hyperbolic metric
 and the conformal structure
on the fibers, respectively. This is reminiscent, as
the following analogies will confirm,
 to
a Kleinian group
acting on the Riemann
sphere.

\subsection{Limit sets}

\begin{definition}
 Consider $s: M \to
{\bf H}M$  a continuous
section (this always exists since ${\bf H}^d$
is contractible).
Let $A$ be a {\bf { closed non compact}}
 subset of $Isom(M)$. The limit set $L_A$ of $A$
is the set of the limits in ${\overline {\bf H}}M$
of the sequences $( \gamma_n(s(x_n)) )$, for
$(x_n)$ a sequence in $M$ and
$( \gamma_n)$ a divergent
sequence of $A$ (i.e. $\{ \gamma_n / n \in {\bf N}\}$
 is a closed discrete subset
of $Isom(M)$).
\end{definition}

\begin{fact}  The definition above
does not depend on
the choice of the  section $s$.
The limit set $L_A$ is a  non empty closed
subset of ${\bf S}^\infty M$.
\end{fact}

\begin{proof} $L_A$ is non-empty since $ \overline {{\bf H}}M$
is compact.
 One easily see that if for a sequence $ (\gamma_n )$,
there is $s_0 \in {\bf H}M$, such that
the elements $\{ \gamma_n(s_0) / n \in {\bf N} \}$ remain
in a compact subset of ${\bf H}M$, then $( \gamma_n )$
is equicontinuous. This shows that $L_A$ is
contained in ${\bf S}^\infty$.
To see that $L_A$
does not depend on the choice of the section $s$, consider
a distance $d$ on ${\bf H}M$
inducing the hyperbolic metric on the fibers. If $s^\prime$
 is another section, then $d(f_ns(x), f_ns^\prime(x))
 = d(s(x), s^\prime(x))$. As usual, this implies that
if $ (\gamma_n(s(x)) )$ converges
in $\overline {\bf H}M$, then,
$ (\gamma_n(s^\prime(x)) )$ also tends to the same limit.
\end{proof}

\begin{proposition} Let $\Gamma$ be a closed non compact
subgroup
 of $Isom(M)$. Then the limit set $L_\Gamma$ is the smallest
$\Gamma$-invariant  closed
subset of ${\bf S}^\infty M$ which  projects
surjectively  onto $M$.
\end{proposition}

\begin{proof} The proof
is the same as
for Kleinian groups.
Let $F $ be  a $\Gamma$-invariant
subset of ${\bf S}^\infty M$ projecting onto $M$.
 One
 naturally
construct
$H(F) \subset
\overline{\bf H}M$,
 the  fiberwise convex hull
 of $F$. We can choose the section $s$ (in the definition
of limit sets) to have image in
$H(F) \cap {\bf H}M$ (since the fibers of $H(F)$ are still
contractible and since $F$
projects onto $M$). Therefore $L_\Gamma \subset H(F)$
since $H(F)$ is closed and invariant. In fact
$L_\Gamma \subset H(F) \cap {\bf S}^\infty M = F$.
\end{proof}

In the definition
of limit sets, we used continuous sections. In fact,
less regular sections may be equally useful.
We will need later the particular following statement

\begin{fact} [measurable sections]  \label{measurable.section}
Suppose that a closed
subgroup
$\Gamma$ of  Isom$(M)$ preserves  a measurable
section $N \to {\bf H}M$, defined over a measurable
$\Gamma$-invariant subset $N$ with a  positive volume. Then,
$\Gamma$  is compact.
\end{fact}

\begin{proof} Let $s: M \to {\bf H}M$ be
a measurable $\Gamma$-invariant section.  Let
$K$ be a compact subset of $M$, along
which $s$ is continuous and such that
 Vol$(K) > 1/2 $ Vol$(M)$.
 Let $(f_n )$ be
sequence
 in $\Gamma$ and  $x_n
\in f_n^{-1}K \cap K  $, this exist
because of the volume condition.
 By continuity of $s \vert K$,
the set $ \{ f_n (s(x_n)) , n \in {\bf N} \}
 \subset {\bf H}M$ is precompact. This
implies that $(f_n )$ is
equicontinuous. Therefore,  $\Gamma$
is compact.
\end{proof}

\subsection{Limit sets and approximately stable
foliations}

\begin{proposition} \label{south.north} Let $(f_n)$ be a sequence
of $Isom(M)$ such that $(f_n^{-1})$
 admits  an approximately stable
foliation ${\cal A S}((f_n^{-1}))$.
 Then the limit
set $L_{ (f_n) }$ is  the image of  a section
$M \to {\bf S}^\infty M$, i.e. a field
of isotropic lines, which is in fact
the normal direction of  ${\cal A S}((f_n^{-1}))$.

Assume furthermore that $(f_n)$ admits
 an approximately stable
foliation. Then if $v \in TM-AS((f_n))$,
$Df_nv $ tends to
$\infty$, and converges projectively (i.e.
after normalization)
to $L_{ (f_n)}$. The convergence is
 uniform
in compact subsets of $TM-AS((f_n))$.

This
convergence is in particular valid in
${\bf H}M$ (since it is missed
by  lightlike hyperplanes
 and hence ${\bf H}M
\subset TM-AS((f_n))$. More precisely,
let $U $  and $V$ be
neighborhoods of respectively
$AS((f_n^{-1} ))^\perp$ and
$AS( (f_n ))^\perp$
in ${\bf S}^\infty$. Then there is
$N$, such for, $f_n( {\bf S}^\infty M -V)
 \subset U$, for $n >N$.
\end{proposition}

\begin{proof}
We argue by contradiction to
prove that
the direction $L_{ (f_n ) }$
is the normal direction of
${\cal A S}((f_n ))$.
 Let $s: M \to {\bf H}M$
be
a section, that is a vector field
$<s(x), s(x)> =-1$ for $x \in M$.
 We have to prove that
$Df_n s$
 converges projectively (i.e. in direction)
 to $AS((f_n^{-1} ))^\perp$.
 If this were  false, it would
exist a sequence $x_n$
such that $D_{x_n}f_n (s(x_n))$
 converges projectively
to a vector $u \notin AS( (f_n^{-1} ))^\perp$.
 Hence,,  $D_{x_n}f_n (s(x_n)) = \alpha_n u_n$,
with $u_n \to u$ and
$\alpha_n \to \infty$ (because
if not $(f_n)$ would be equicontinuous).
Thus $Df_n^{-1}u_n = \alpha_n^{-1}s(x_n)$.
 Thus,  by definition $u$
 is strongly approximately
stable for $(f_n^{-1} )$, and hence by
\ref{cocycle.lorentz}, belongs to $AS( (f_n^{-1} ))^\perp$,
 which yields a
 contradiction.
The last   statement may be checked in a similar
fashion.
\end{proof}

From the above theorem, one can conclude  the two following
``continuity'' corollaries:

\begin{corollary} \label{continuity.field} Let $(f_n)$ be a sequence
of $Isom(M)$ such that $(f_n)$
 and
 $(f_n^{-1})$
 admit  approximately stable
foliations  ${\cal A S}( (f_n))$ and ${\cal A S}((f_n^{-1}))$.
Let $X$ (resp. $Y$) be  a continuous  vector
field tangent to $AS((f_n))^\perp$ (resp.
$AS((f_n)$). Then there is sequence $(X_n)$
 (resp. $(Y_n)$
of continuous vector fields converging
in the $C^0$ topology  to $X$
(resp. $Y$), and such that $(Df_nX_n)$
 converges in the $C^0$ topology to $0$ and
$(Df_nY_n)$ is bounded.

\end{corollary}

\begin{proof}  We will only prove the claim for $X$,
as the same argument works for $Y$.
  To simplify notation,
 we will do the proof for
$g_n = f_n^{-1}$ instead of $(f_n)$.  By the Proposition
 above, after passing
to subsequences,  we may find   sequences of
 neighborhoods  $(U_n)$
 and  $(V_n)$,
collapsing to $AS((f_n))^\perp$ and
$AS((f_n^{-1}))^\perp$,  respectively,
 and such that $Df_n(TM-U_n) \subset V_n$. Now, let $(X_n^\prime)$
 be a sequence of vector fields such that $X_n^\prime \notin U_n$, but
$(X_n^\prime)$ converges to some
non singular vector field tangent to $AS((f_n))^\perp$. Thus,
from the Proposition above, we infer that  $X_n^{\prime \prime}=
Df_nX_n^\prime$
is a vector
field with a big norm and converging in direction
in $AS((f_n^{-1}))^\perp$.  Therefore, for suitable
distortion functions $\lambda_n$, going uniformly to $\infty$,
we  may get  $X_n =
\lambda_n X_n^{\prime \prime} \to X$. But
$Dg_nX_n = (1/\lambda_n)X_n^\prime \to 0$.
\end{proof}

\begin{corollary} \label{continuity.value} Let $A$ be a closed
non compact subset of Isom$M$ and $u \in L_A \cap
 {\bf S}_x^\infty M$.
 Then,  there is a continuous (in fact Lipschitz) section
$\sigma : M \to L_A$ with
$\sigma(x) = u$.
\end{corollary}

\begin{proof} $u$ is a  (projective) limit
of
a sequence $ (f_n (s(x_n)) )$
where $(f_n)$ is a sequence in $A$
 and $x_n \to x$ (and $s$ is fixed
section of  ${\bf H}M$).  After
passing
to a subsequence we may
 assume that $( f_n^{-1} )$
admits an approximately
stable foliation. We set
$\sigma (x)
=  AS((f_n^{-1}))^\perp$.
 From the Proposition above, we conclude that $\sigma(x) = u$.
\end{proof}

\subsection{End of the proof of Theorem
\protect\ref{theo.generalized}}

\label{last.part}

  The
 statement in Theorem \ref{theo.generalized},
about  the dynamics of $(D \phi_n)$ on $TM$,
  follows from
Proposition \ref{south.north}.

To finish the proof of Theorem \ref{theo.isometry}, it
 remains to check its  ergodicity
statement.
 For this end,  to simplify notation,
suppose that $\phi_n = f_n$.
let $\sigma: O \to {\bf R}$ be
a continuous  function, invariant by all
the $f_n$.

	Consider a vector field $X_n$ as in \ref{continuity.field}. Let
$\alpha$ be a real such that $\exp_x \alpha X_n(x)$ exist
for any $x \in M$.
Then,  $\sigma(\exp_x\alpha X_n(x))
-\sigma(x)= \sigma( \exp_{f_nx}\alpha Df_nX_n(x))
-\sigma(f_nx)$.  Denote this difference by
$\Delta_n(x)$.

Let $(x_n)$ be sequence
such that
$f_nx_n$ belong
to the same compact $K$. Then
$\Delta_n(x_n) \to 0$,
as
$Df_nX_n \to 0$, in the $C^0$ topology and because
$\sigma$ is uniformly continuous on $K$.

Suppose that $x_n \to x$ and
 $X_n(x_n) \to u \in AS(f)^\perp$. Then
$\sigma(\exp_x \alpha u) -\sigma(x) =
\lim (\sigma(exp_{x_n} \alpha u) -\sigma(x_n)) = \lim \Delta_n(x_n)
= 0$,
 for any $\alpha$ such that $\exp_x \alpha u \in O$.
 That is $\sigma$
 is locally constant
along the leaf of $AS(f)^\perp$ passing through
$x$. Remember finally that, from \ref{def.non.escaping},
almost every $x$ is a limit
of a sequence such as $(x_n)$,
for some compact $K$.

\section{Foliation compactification} \label{foliation.compactification}

Here we define a topology
  on $\hbox{Isom}(M) \cup
{\cal F G}$ and then an
ideal boundary
$\partial_\infty \Gamma$ for $\Gamma$
a closed non compact subgroup
of $\hbox{Isom}(M)$.

We fix a section $s: M \to {\bf H}M$
and we choose
a distance $d$ on $\overline{\bf H}M$ (in
fact the choice
of $s$ permit to construct a natural Euclidean
 fiberwise distance, we then
tensorize by
a distance on $M$ and thus  we
 get a semi-canonical distance
 $d$). Now
we embed  Isom$M \cup {\cal F G}$
 in the  Sect$\overline{\bf H}M$, the
space of sections of $\overline{\bf H}M$.
To $f \in $ Isom$M$, we associate
the section $fs$, and to a foliation
${\cal F}$
 we associate its normal direction
field ${\cal F}^\perp$. We then endow
Isom$M \cup {\cal F G}$ with the $C^0$
topology defined by the
distance $\rho( \sigma_1, \sigma_2) =
 \sup_{x \in M} d(\sigma_1(x), \sigma_2(x) )$.

From \ref{south.north}, we deduce:

\begin{fact} If a sequence $(f_n)$ of
Isom$M$
 converges to ${\cal F} \in
{\cal F G }$, then
${\cal F} = {\cal A S}((f_n^{-1} ))$.
\end{fact}

Since ${\cal F G}$
is compact, we have:

\begin{corollary} Isom$M \cup {\cal F G}$
is a (metrizable) compact space endowed
with a natural continuous action
of  Isom$M$.

\end{corollary}

\begin{definition} Let $\Gamma$ be a closed non
compact subgroup of
 Isom$M$.  Its boundary $\partial_\infty \Gamma$
 is the intersection with ${\cal F G}$
 of the closure of $\Gamma$
in Isom$M \cup {\cal F G}$.

\end{definition}




\section{Preliminaries  on elementary groups. Proof of Theorem
\protect\ref{theo.amenable} } \label{elementary}

Let
${\cal F}^1 \cup {\cal F}^2$ be
a bi-foliation  (Definition \ref{def.bifoliation}).
Its  set of tangent foliations
 Apa$({\cal F}^1, {\cal F}^2)$
has the following description. Let $C$ be the coincidence set
$\{x \in M / T_x{\cal F}^1 = T_x {\cal F}^2 \}$.
Let $\pi_0(M-C)$ be the set of connected components of
$M-C$, and consider a map: $c: \pi_0(M-C) \to \{1, 2 \}$. This
allows to construct a foliation, ${\cal F}_c$, by the rule:
$T_x{\cal F}_c = T_x {\cal F}^{c(U)}$, in the connected component
$U$ of $M-C$, and $T_x{\cal F}_c = T_x{\cal F}^1 = T_x {\cal F}^2$
in $C$.  It is straightforward to
verify that ${\cal F}_c \in $ Apa$({\cal F}^1, {\cal F}^2)$
(in particular there is no loss of continuity or
Lipschitz character in this construction). In fact
we have a topological identification
 of Apa$({\cal F}^1, {\cal F}^2)$ with
$\{1, 2\}^{ \pi_0(M-C)}$.

Observe that we have a natural distance in Apa$({\cal F}^1, {\cal  F}^2)$
 defined by $d({\cal F}_c, {\cal F}_{c^\prime})
= \Sigma_{U \in  \pi_0(M-C)} (c(U)-c^\prime(U)) \hbox{Vol}(U)$. Any group
preserving the volume and the bi-foliation, preserves this distance.

\begin{fact}
Let $\Gamma \subset $ Isom$(M)$ be
a closed non compact subgroup
with the property that
for any $x \in M$,
 $L_\Gamma(x)$
 $= L_\Gamma \cap {\overline {\bf H}}M = {\bf S}^\infty M$,
has exactly  1
or 2 points. Then,  $\Gamma$ is elementary.
\end{fact}

\begin{proof} Let $u \in L_\Gamma(x)$. From the corollary above
there is a foliation ${\cal F}(u)$, such that the normal
direction ${\cal F}^\perp$ is contained in
$L_\Gamma$ and  equals
$u$ at $x$. Let $C$ be
the coincidence set of all these foliations. Obviously, it
equals the set of points $x$  where $L_\Gamma(x)$ has cardinality 1.
Let $U$ be a component of $M-C$, then the restrictions to
$U$ of two foliations ${\cal F}(u)$
and ${\cal F}(u^\prime)$ are identical or everywhere transverse.
Therefore, as above, one may construct two foliations
${\cal F}^1$ and ${\cal F}^2$, which are transverse in each component
$U$, and each of them equals a foliation of the form
${\cal F}(u)$ (in the component $U$). Therefore,
$L_\Gamma$ can be
defined by a bi-foliation.
\end{proof}

\paragraph{Proof of Theorem \protect\ref{theo.amenable}}
Let $\Gamma \subset $ Isom$(M)$ be
a closed non compact amenable subgroup. To
prove that it is elementary, we will check that
it satisfies the condition of the fact above.
 By amenability, there
is a $\Gamma$-invariant probability
measure $\mu$ onto  ${\bf S}^\infty M$, projecting
 onto the Lorentz measure
of $M$. Let $\mu_x$ the conditionals
of $\mu$ on the fibers ${\bf S}^\infty_x M$.

Let us show that for
almost all $x \in M$ (in the sense of the Lorentz
measure), the support
of $\mu_x$ has exactly
1 or 2 points. Indeed, if not a fiberwise  barycenter
construction (see \cite{} for the classical one)
yields a measurable $\Gamma$-invariant
section $\sigma:N \subset  M \to {\bf H}M$
over a $\Gamma$
 invariant subset $N$. This is impossible
by Fact \ref{measurable.section}.
$\Box$

\paragraph{The bi-foliation of an
elementary group.}

Observe that an elementary group may preserve many
 bi-foliations.  However
we
have the following fact allowing to define the bi-foliation associated
to an elementary group:

\begin{fact} Let $\Gamma$ be an elementary group. Then
$\Gamma$
preserves a unique  bi-foliation ${\cal F}^1 \cup {\cal F}^2$
 determining a maximal apartment among that determined by all
the bi-foliations preserved by $\Gamma$. It is characterized by:
$L_\Gamma(x) = \{ (T_x{\cal F}^1)^\perp, (T_x {\cal F}^2)^\perp \}$
for any $x \in M$. We call it {\bf the bi-foliation of} $\Gamma$.

\end{fact}

\begin{proof} For a bi-foliation
${\cal G}^1 \cup {\cal G}^2$, denote by
 $T( {\cal G}^1, {\cal G}^2) = \{x \in M/ T_x {\cal G}^1
\neq T_x{\cal G}^2 \}$ its transversality (open) set. Now, we
will consider only bi-foliations preserved by
$\Gamma$.  Let ${\cal G}^1 \cup {\cal G}^2$
and ${\cal L}^1 \cup {\cal L}^2$ be two of them. Observe
then that along
$T({\cal G}^1, {\cal G}^2)$, ${\cal L}^1$
and
${\cal L}^2$
are tangent to ${\cal G}^1 \cup {\cal G}^2$. Indeed, the opposite
situation, would give a $\Gamma$-invariant  $3-$ or
$4-$valued section of ${\bf S}^\infty M$
along some open subset of $T({\cal G}^1,  {\cal G}^2)$. As above,
the barycenter construction (the classical one in this case),
yields a $\Gamma$-invariant section of ${\bf H}M$, over
some open subset of $M$. This implies
$\Gamma$
is compact, by  Fact \ref{measurable.section}.

Observe now that by the above gluing process, one may
construct
a bi-foliation ${\cal M}^1 \cup {\cal M}^2$
 such that $T({\cal M}^1 ,  {\cal M}^2)
=T({\cal G}^1,{\cal G}^2) \cup T({\cal L}^1, {\cal L}^2)$. We
construct the wished maximal foliation  as one having
a maximal
transversality set. This exists by
compactness of
the space of bi-foliations.
\end{proof}


\begin{remark}
{ \em
 We close  this section by   making  some remarks on the
 general set-up of the ``bundle compactification''
(\S \ref{bundle.compactification}), when
we let a generalized dynamical
system
of $C^1$ diffeomorphisms acting
 on the bundle
$\overline {S_k} M \to M$.
 One
considers a given measure $\nu$
preserved by
the dynamical system, and choose
an invariant  measure $\mu$ on
 $\overline {S_k} M$
of maximal support  and  projecting onto  $\nu$.
 It may happen (in fact generically)
that $\nu$
has a full measure in
 $\overline {S_k }M$, or even
worse, $\nu$ can be ergodic. If not, i.e.
  $\nu$ has support inside  ${\bf S}_k M$
 (that is the support of $\nu$ does not
interest $\partial_\infty {\overline {\bf S}}_k M$)
then,  the barycenter construction
(in general
 Hadamard
 spaces) yields an invariant
measurable metric. However, if
$ k >2$, the universal
 symmetric space $S_k =
 SL(k, {\bf R})/SO(k)$ is not of
 negative curvature
(i.e. it has
a higher rank), and there is
 no way to construct
barycenters for  measures
 supported in the Hadamard
 boundary. In fact, there
 are alternative boundaries
which  may be efficient
in this matter. Depending
 on the interpretation
of such a boundary, that is by
modeling it as a kind of
flag spaces,
the construction of
limit sets
 yields, roughly speaking, flag-fields
on  $M$.
 For example, the oseldec's decomposition for
a diffeomorphism $f$, may be handled by
looking
to points in the limit set (of
thegroup $\{ f^n, n \in {\bf Z}  \}$),
with special ``approach''.  For instance,
the conical and horospherical
limit
sets may be interesting in
this regard (see \cite{Sul}).
}
\end{remark}

\section{Partial hyperbolicity} \label{partial.hyperbolicity}

\begin{definition} \label{def.hyperbolic}
An  elementary
group
$\Gamma$ is called {\bf partially hyperbolic}
if
its  associated bi-foliation
${\cal F}^1 \cup  {\cal F}^2$ is non trivial, that is
 the (open) transversality
locus $T= \{x \in M \, / \, T_x{\cal F}^1 \neq T_x{\cal F}^2 \}$
is non-empty,
or equivalently, the (closed) coincidence locus
$C =  \{x \in M \, / \, T_x{\cal F}^1 = T_x{\cal F}^2 \}$
 is a proper subset of $M$.
 \end{definition}

Now, we will justify the word ``partially hyperbolic'' in a
 dynamical viewpoint
(of course, partially hyperbolic,  is also reminiscent to
 the term ``elementary
hyperbolic'' in the theory of Fuschian groups). We keep the
 notations
above. Denote by ${\cal N}^1$ and
${\cal N}^2$  the two (one dimensional) normal
directions of ${\cal F}^1$ and ${\cal F}^2$,  respectively. We may
assume that they are
 orientable,  after passing to a finite covering. Let
$X^1$ and $X^2$ be two Lipschitz non singular vector-fields
 orienting ${\cal N}^1$ and ${\cal N}^2$,
respectively.

The group $\Gamma$ respects the set of  (open)
connected components of the transversality
locus $T$. Since these components have positive Lorentz
volume and since $\Gamma$ preserves the volume,  then
each component $U$ is preserved by a finite index subgroup
of $\Gamma$,   say $\Gamma$ itself. It then follows
that $\Gamma$ preserves the closure $\bar{U}$ and
the directions of $X^1$ and $X^2$, along it. Therefore we get
two derivative cocycles: $\lambda^1$ and
$\lambda^2: \bar{U} \times \Gamma \to {\bf R}$, defined
by $D_xf(X^i(x)) = \lambda^i(f, x)X^i(f(x))$, for
$i \in \{1, 2\}$.
This section is devoted to the proof of:

\begin{theorem} \label{theo.hyperbolic}
Let $\Gamma$ be a partially hyperbolic elementary group.
Then,  up to a subgroup of finite index, $\Gamma$ is a direct product of
a compact group by ${\bf Z}$ or ${\bf R}$.
Furthermore, we can find two foliations
${\cal L}^1$ and ${\cal L}^2$ generating the
same bi-foliation ${\cal F}^1 \cup {\cal F}^2$, and defining two cocycles
$c^1$ and $c^2$ satisfying
the following condition.
Let $g$
be a non trivial element of the ${\bf Z}$ or the ${\bf R}$
part of $\Gamma$
 and $K$ a compact subset in the transversality set of the
bi-foliation. Then,  there is an integer
 $p = p(K) >0$,  such that for
 $f = g^p$,   we have:

$c^1(f^n, x) <1/2$, if  $n >0$, $x \in K$ and $f^nx \in K$, and

$c^2(f^n,x)>2$, if  $n > 0$, $x \in K$ and  $f^nx \in K$.

Finally, in the transversality set $T$,
the normal foliation of ${\cal L}^1$
and ${\cal L}^2$
 are  respectively the negative and positive Liapunov spaces
of $g$. In particular,  ${\cal L}^1$
and ${\cal L}^2$ are
preserved by
$g$.

\end{theorem}

\paragraph{Beginning of the proof.}
Consider
an auxiliary complete distance $\rho$ on $\Gamma$. For example
embed $\Gamma$
 is the frame bundle $P_M$, and  take $\rho$ to be the restriction
of a  distance on $P_M$ induced by a complete Riemannian metric. We denote
open balls around the identity by $B^\rho(1,r)$.

\begin{fact}  \label{fact.punctual}
 With the
notations above let
$K$ be a compact of $U$.  There are:

(i) a function: $r \in {\bf R}^+ \to c_r^K \in {\bf R}^+ $  such
that $\lim_{r \to \infty} c_r^K = \infty$, and

(ii) a map $s: K \times \Gamma \to \{1, 2\}$, defined
for $(f, x)$   such that $f(x) \in K$. We denote
$u(f, x) = s(f, x)+1 $ mod 2.

These maps satisfy the condition
that whenever,
$ f \notin B^\rho(1, r)$,  then:

$$ \lambda^{s(f,x)}(f, x) < (c_r^K)^{-1}
\mbox{ and} \;\; \lambda^{u(f,x)}(f, x) >c_r^K$$
\end{fact}

\begin{proof}
In $U$, $\Gamma$ preserves the directions
${\cal N}^1$ and ${\cal N}^2$, and hence also the orthogonal
$({\cal N}^1 \bigoplus {\cal N}^2) ^\perp$. This last space
is spacelike, i.e. the Lorentz metric restricted to it is
positive definite. Therefore, in the  compact $K$ of $U$,
we have an uniform bound
of the restriction  $D_xf \vert ({\cal N}^1 \bigoplus {\cal N}^2) ^\perp$
 and $(D_xf)^{-1} \vert ({\cal N}^1 \bigoplus {\cal N}^2) ^\perp$,
for $x \in K$ and $f \in \Gamma$, such that $f(x) \in K$.

It then follows by the volume preservation property that
the product
$\lambda^1(f, x) \lambda^2(f, x)$ belongs to some fixed compact
interval (around 1)
in $]0, \infty [$, for $x \in K$, and $f \in \Gamma$, with
$f(x) \in K$. That is the cocycles $\lambda^1$ and
$\lambda^2$ are almost one the inverse of
the other, provided we restrict
ourself
to $K$.

 To prove the
estimates contained in the statement of the
  Fact,
  we argue by
 contradiction.
Suppose that for a divergent sequence  $(f_n)$
in $\Gamma$, there is a sequence $(x_n)$ of
points of $K$, such that $f_n(x_n) \in K$ and such
 both
$\lambda^1(f_n, x_n)$ and $\lambda^2(f_n, x_n)$
remain bounded.
Hence,  $(D_{x_n}f_n )$
 is equicontinuous and thus it follows
from \ref{start.unimodular},that  $(f_n)$ is
not divergent !
Therefore, there is a function $c_r^K$,
satisfying the
property (i) of the Fact,
 and  such that for any $x$
 and $f$ as in the Fact, there is  some $
s = s(f, x) \in
\{1, 2\}$, such that the
following inequalities hold:  $\lambda^s(f, x) < (c_r^K)^{-1}$.
and
 $\lambda^u(f, x) >c_r^K$
 with $u = s+1$ mod 2.
\end{proof}

\begin{fact} \label{fact.north.south}

There is a partition $f = A_1 \cup A_2$ satisfying
the following conditions.
Let $U_1$ (resp. $U_2$) be a neighborhood
in (the projective isotropic cone) ${\bf S}^\infty M$
of ${\cal  N}^1$ (resp. ${\cal  N}^2$)
 along the closure $\bar{U}$.

There is $r$ such for $f \in \Gamma- B^\rho(1,r)$, we have:
$$ \hbox{if} \;  f \in A_1,  \;
\hbox{then}:  Df( {\bf S}^\infty M-U_1) \subset U_2, \; \hbox{and}, $$
$$ \hbox{if} \; f \in A_2,  \;
\hbox{then}:  Df( {\bf S}^\infty M-U_2) \subset U_1.$$

\end{fact}

\begin{proof}
It is easily seen that one
 can localize the  compactification of  isometry
groups of $M$ to that of $\bar{U}$. So, here the boundary
of $\Gamma$ (acting on $\bar{U}$) consists
of the two foliations ${\cal F}^1 \vert \bar{U}$
 and ${\cal F}^2 \vert  \bar{U}$.  Let
$B_1$ and $B_2$ be two disjoint neighborhoods
of these last foliations in the compact $ \overline{ \hbox{Isom}
(\bar{U}) }$. So, for $r$ big enough, every element
$f$, with $ \rho(1, f) >r$, belongs to exactly
one of the neighborhoods $B_1$ or $B_2$.
Suppose that we can't find $r$ (big enough) satisfying
the claim. Then there is a divergent sequence $(f_n)$
 contained in $B_1$ or $B_2$ (say $B_1$) and
 not satisfying the desired inclusions. Thus the approximately
stable foliation
of $(f_n)$ is ${\cal F}^1 \vert {\bar{U}}$. Now we apply
the expanding
property of the approximately stable foliation \ref{south.north}. This
would lead to a contradiction if we check that for some
subsequence $( \phi_n )$ of $(f_n)$, the approximately
stable foliation ${\cal A S }( \phi_n^{-1})$ is
${\cal F}^2$.   For this last fact,
 we just remark
that the opposite situation, i.e. ${\cal A S }((f_n^{-1})) = {\cal F}^1$
is impossible, indeed this would imply that
 the complementary of
a small neighborhood of ${\cal N}^1$ is mapped by
$Df_n$ to a small neighborhood of ${\cal N}^2$. This contradicts
the fact that $Df_n$ preserves ${\cal N}^2$.

Finally,  extend arbitrarily the partial partition
$B_1 \cup B_2$ to a partition $\Gamma = A_1 \cup A_2$.
 \end{proof}

\begin{fact} $s(f, x)$ and $u(f, x)$ are independent
of $x$.

\end{fact}

\begin{proof}  It is clear that  by choosing $U_1$ and $U_2$ small
enough that for  $f \in A_2$  far from
the identity map, $s(f, x)$ cannot
  be $2$ for
any $x \in K$. Therefore,  $s(f) =1$,  and thus
 does not depend on $x$.
\end{proof}

\paragraph{Liapunov exponents.}
Henceforth, we will assume that  the following hypothesis holds: \\

{\bf Hypothesis}.  Vol$(K) >$ (1/2)Vol$(U)$
and
 choose $r_0$ such that $c_{r_0}^K >2$.

\begin{lemma} \label{lemma.semi.stability}
Let  $f$ and $g$ be two elements of $\Gamma$ outside
the ball $B^\rho(1,r_0)$ such that $ s(f) = s(g)$. Then there is $x \in K$
such that $f(x)$ and $gf(x)$ belong to $K$
and $\lambda^{s(f)} (gf, x) <1/4$.
 In particular, if  $gf \notin B^\rho(1, r_0)$, then $s(gf) = s(f)
$ ($= s(g)$). Furthermore,  if $f^n \notin B^\rho(1, r_0)$ for
all $n>0$, then $s(f^n) = s(f)$ for all $n>0$.
\end{lemma}

\begin{proof} Observe that if three subsets $A$,
$B$ and $C$ of $U$ have volume $> $ (1/2)Vol$(U)$, then
$A \cap B \cap C \neq \emptyset$.
Apply this to $K$, $f^{-1}K$
 and $gf^{-1}K$, we get a point $x \in K$ such that
$f(x) \in K$, and $gf(x) \in K$. Thus
 $\lambda^1(gf,x) = \lambda^1(g,f(x)) \lambda^1(f,x)$, and
hence $\lambda^1(gf,x) <1/4$, if $f$ and $g$
satisfy the conditions of the Lemma. Of course
if $gf \notin B^\rho(1, r_0)$, then $s(gf) \neq 2$,
and thus $s(gf) =1$.
 Using this, the last part of the lemma
is  proved by  induction.
\end{proof}

Consider $l^1(f, x) = \lim_{n \to + \infty}
(\log \lambda^1(f^n,x))/n$.
 Define analogously $l^2(f, x)$. They are
Liapunov exponents and thus  exist
 almost everywhere.

\begin{fact}  \label{fact.semi.stability}
 Let $f$ be such that $f^n \notin B^\rho(1, r_0)$ for
all $n>0$. Suppose that  $s(f) =1$. Then,
for almost every $x \in U$,
$\l^1(f, x) <0$, that is ${\cal N}^1$ is
 the negative Liapunov space of $f \vert U$. Furthermore,
 $\int_U l^1(f, x)dx <(-\log2)$Vol$(K)$.

\end{fact}

\begin{proof}
To $x \in K$, associate its sequence of
positive return times $(n_i(x) )_{i \in {\bf N}}$.
To simplify notation, fix $x$ and denote the sequence
by $(n_i)$. Thus, $f^{n_i}x \in K$. From the Lemma
 \ref{lemma.semi.stability},
$ s(f^n)=1$
for all $n >0$, and hence $\lambda^1(f^{n_{i+1}-n_i}, f^{n_i}x)
<1/2$. Thus,  by the cocycle property of $\lambda^1$,
$\lambda^1(f^{n_i}, x) <(1/2)^i$, and
hence,   $\log \lambda^1(f^{n_i}, x) / n_i <i/n_i (-\log2)$.
Let $\chi_K$ denote
 the characteristic function of $K$.  Observe that
$i/n_i$ equals the partial Birkhoff sum
$(\chi_K(x) + \chi_K(fx)+ \ldots \chi_K(f^{n_i}x) )/n_i$.
 So, we have proved
$l^1(f,x) <(-\log2) \chi_K^*(x)$, where
 $\chi_K^*$  stands for
 the Birkhoff sum of $\chi_K$.  In
particular $l^1(f, x)<0$ for almost every $x \in K$.
Let $K^*$ be the saturation of $K$, by $f$,
$l^1(f,x) < (-\log2) \chi_K^*(x)$, for $x \in K^*$, because
both of the two functions $l^1(f,x)$ and $\chi_K^*$
are $f$-invariant.

In particular: $\int_{K^*} l^1(f,x)dx (-\log2) \int_{K^*} \chi^*
 = \int_U \chi_K^* = $ Vol$(K)$.

It remains to  prove that $l^1(f,x)<0$, almost everywhere in $U$.
To this end, observe that
the sequence generated by $f$ is divergent because $f$
has non vanishing exponents. Hence if we replace
$K$ by a bigger compact $K^\prime$, then for some
 $g= f^p, p>0$, the powers $\{g^n, n>0\}$ lie outside
the ball analogous to $B^\rho(1, r_0)$
associated to $K^\prime$. Obviously, the index
$s^{K^\prime}(g)$ is the same as $s(f)$ and
thus equals $1$. Therefore, almost
everywhere in
$K^\prime$, $l^1(f^p,x) <0$. But $l^1(f,x)
= l^1(f^p,x)/p$.
\end{proof}

Consider now the map $ \Lambda^1: \Gamma \to  {\bf R}$,
$\Lambda^1(f) = \int_U \log \lambda^1(f, x)dx$.

\begin{fact} \label{homomorphism}
$\Lambda^1$ is a homomorphism and
 for $f$ such that $f^n \notin B^\rho(1, r_0)$ for
all $n>0$, we have
 $\vert \Lambda^1(f) \vert > \log2$ Vol$(K)$.  Furthermore,
$\Lambda^1 (f) \neq 0$  if and only if
$f$ generates  a non  precompact group (that
is  $ \overline{   \{f^n / n \in {\bf Z} \}  }$
is non compact).
\end{fact}

\begin{proof} We have $\log \lambda^1(fg, x) = \log \lambda^1(f, gx)
+ \log \lambda^1(f,x)$. Thus $\Lambda^1$ is a homomorphism
because $g$ is volume preserving and hence: $\int_U \log\lambda^1(f, gx)
= \int_U \log \lambda^1(f,x)$.

The remaining parts of the fact are obvious or
follow from the preceding Fact.
\end{proof}

\begin{corollary} \label{when.compact}
 Suppose that $Ker\Lambda^1$ is
compact. Then,  up to a subgroup
of finite index, $\Gamma$
is a direct product of a compact
group by ${\bf Z}$ or ${\bf R}$

\end{corollary}

\begin{proof} Suppose to start with that
 $Ker \Lambda^1$ is trivial, that
is every element $f \in \Gamma$ generates a non precompact group.
 Therefore, from  Fact \ref{homomorphism}
 $\Lambda^1$ is injective, and thus $\Gamma$
is abelian and torsion free. Let  $G$ be a closed subgroup
of
$\Gamma$ isomorphic to ${\bf Z}^2$.  It is obvious that
for a  big enough radius $r$,  $f \in
G-B^\rho(1, r)  \Longrightarrow f^n \notin B^\rho(1, r_0)$. But,
from Lemma \ref{lemma.semi.stability},
 $\vert \Lambda \vert$ is bounded from below for elements
satisfying this condition. This means that $\Lambda^1: G
\to {\bf R}$ is proper, which is impossible. Thus,
because it can't  contain  a closed copy of ${\bf Z}^2$,
$\Gamma$
 must  be isomorphic to
${\bf Z}$ or ${\bf R}$ (remember that it is torsion free).

Now when  $Ker\Lambda^1$ is merely compact, we
 may argue with the quotient $\Gamma / Ker\Lambda^1$
 which enjoys
the same properties
as $\Gamma$. So, we obtain
that $\Gamma$ is a semi-direct product of
${\bf Z}$ or ${\bf R}$ by a (normal) compact
group.
We may find in a standard way
 a subgroup of finite index which is
a direct product of
 ${\bf Z}$
or ${\bf R}$
  by a compact group.
\end{proof}

\paragraph{Torsion.}
We  now check  that $Ker \Lambda^1$
is compact. Equivalently, we suppose that $\Lambda^1 =0$
and thus show that $\Gamma$ is compact. This is   standard
for the identity  component of $\Gamma$. Indeed,  a Lie group
for which every element generates a precompact subgroup,
is compact (see for instance \cite{D'A}).
 Hence, without loss of generality,
we may
restrict ourselves to the case where $\Gamma$ is discrete.
 It is thus a torsion group, and we have to check it is finite.

Let $f$ be an element of order
$k$. Consider
$I_1(f) = \{i \in \{1, \ldots, k-1\} / f^i \notin B^\rho(1, r_0),
s(f) =1 \}$. Define analogously $I_2(f)$.

From  Lemma \ref{lemma.semi.stability},
 we have the following ``semi-group'' property:

\begin{fact} \label{semi.group}
 Let $\alpha, \beta \in I_1(f)$, then
$\alpha+ \beta \in I_1(f) $, unless,
$f^{\alpha+ \beta} \in B^\rho(1, r_0)$. The same statement holds
for
$I_2(f)$.

\end{fact}

Consider $P = \{f^n/ f \in B^\rho(1, r_0),
n \in {\bf Z} \}$. Since
$B^\rho(1, r_0)$ is finite and any element
has finite order, it  follows that $P$ is finite.

We deduce from the Fact above  that
every (finite cyclic) group
 intersects  non trivially the ball $B^\rho(1, r_0)$.  Indeed,
if not, we would obtain a partition $\{1, \ldots, k-1\} = I_1(f)
\cup I_2(f)$, where $f$ is a generator of order $k$ of the given
cyclic subgroup.  The Fact above implies that only one
part, say $I_1(f)$ is non empty. So, we apply,   Lemma
 \ref{lemma.semi.stability}
to $f$ and $f^{k-1}$ and we obtain that $\lambda^1( \hbox{Identity}, x)
<1/4$ which is obviously impossible.

If  $f$ is of prime order, then, all its
non  trivial powers generate the same group, and
hence $f$ is a power of some element of
 $B^\rho(1, r_0)$, that
is $f \in P$.

To treat the general case,
represent the congruence group ${\bf Z}/k{\bf Z}$
as $\{\dot{0},  \dot{1}, \ldots, \dot{(k-1)} \}$. Let
${\cal G} =
\{ \dot{\alpha} / (\alpha,k) = 1 \}$ be the set of generators of
 ${\bf Z}/k{\bf Z}$.

Let $f$ be an element of order $k$ that does not belong to $P$,
then  $f^\alpha \notin B^\rho(1, r_0)$, for
$\dot{\alpha} \in {\cal G}$. Hence,  we have
a partition ${\cal G}  = ({\cal G} \cap I_1(f)) \cup
 ({\cal G}\cap I_2(f))$, as above.

Suppose that  $k= p^m$ with
$p$ prime and $\neq 2$.
   Then   $\dot{ \alpha}  \in {\cal G}$ if and only
$\alpha$ is not a multiple of $p$.
Therefore,  we can't
have
$\dot{\alpha}$ and $\dot{\alpha+1}$ in ${\cal G}$,
 noting that  $p \neq 2$. Also, suppose that
$\dot{1} \in I_1(f)$, then   $\dot{2} \in I_2(f)$. Let
$\alpha$ be the smallest number  such that $\dot{\alpha} + \dot{1}
\notin I_1(f)$. Then necessarily, $\dot{\alpha}
 \in I_1(f)$ and  $\dot{ \alpha} + \dot{1
} \notin {\cal G}$, and hence $\dot{\alpha} + \dot{2} \in
I_1(f)$. By recurrence,  we get
$I_2(f) \cap {\cal G} = \emptyset$.  In particular $f$ and
$f^{k-1}$ belong to $I_1(f)$.
We get a contradiction
as  inthe above  case where $k$ was  prime.

To treat the case $p=2$, the previous combinatorial
approach
fails.  However, this may be adapted, if we suppose
a stronger `` semi-group'' property of $I_1(f)$
and $I_2(f)$, involving three instead of two elements. That
is given three elements $\dot{\alpha}$
$\dot{\beta}$ and $\dot{\gamma}$ of $I_1(f)$, then each of
the elements
$\dot{ \alpha} + \dot{\beta}$, and $\dot{\alpha}+ \dot{\beta}
+ \dot{\gamma}$ belongs to $I_1(f)$, unless it belongs
to $B^\rho(1, r_0)$. To have this,
we choose $K$
  in the Lemma
 \ref{lemma.semi.stability},
with
  relative big volume (that is  $\hbox{Vol}(K) >3/4 \hbox{Vol}(U)$), and
 we obtain a statement
(of the Lemma) involving three elements of $\Gamma$.

We have thus proved that $P$ contains all
the elements  having  order of the form $p^m$
with $p$ prime.  One  may push forward the combinatorial argument to
prove that $P = \Gamma$. Instead, we prefer to argue as follows.
Since every cyclic group is generated by groups
with order of the form $p^m$, we deduce that every
subgroup $\Gamma^\prime$ of $\Gamma$ is generated by its intersection
with  $P$.  In particular,
$\Gamma^\prime$ is a finitely generated torsion group. Consider
the adjoint action of $\Gamma$ on itself. It preserves the finite
set $P$. Its Kernel $\Gamma^\prime$ is of finite index in $\Gamma$,
and centralizes $\Gamma$ because it centralizes
the generating set $P$.
 Therefore,  $\Gamma^\prime$ is finite as it is a finitely
generated  abelian torsion group. It then follows that $\Gamma$
is finite. $\Box$

\paragraph{End of the proof of Theorem \protect\ref{theo.hyperbolic}.}
Let $g \in \Gamma$  be as in  Theorem \ref{theo.hyperbolic}.
From the previous development, for each component
$U$ of the transversality set $T$, we can associate
$s(U) \in \{1, 2\}$, such that the normal
direction of
${\cal F}^{s(U)} \vert U$ is the negative Liapunov
space of the restriction $g \vert U$. As in \S \ref{elementary}, this
allows us  to construct two foliations ${\cal L}^1$
 and ${\cal L}^2$, elements of Apa$({\cal F}^1, {\cal F}^2)$
 whose normal directions are  the negative
and positive spaces of $g \vert T$, respectively. Now let $K$ be
a compact subset of $T$, then it meets only finitely
many components $U$, and hence for a sufficiently big
positive integer $p$, we have the estimates stated in  the
theorem for the power $f = g^p$. This ends  the
proof of the theorem.

\section{Proofs of Theorems \protect\ref{theo.elementary},
\protect\ref{theo.isometry} and \protect\ref{theo.flow} }

\paragraph{The case where the bi-foliation of $\Gamma$ is trivial.}

This means there
is a foliation ${\cal F}^1$ such that everywhere $L_\Gamma(x) =
 T_x{\cal F}^\perp$.  It then follows from
 Proposition \ref{south.north}, that the approximately stable
 foliation of any sequence of of $\Gamma$
 is ${\cal F}^1$ and therefore the boundary of $\Gamma$ consists exactly
of ${\cal F}^1$.

To chow the vanishing of the entropy of the elements of $\Gamma$, we
argue by
contradiction. Suppose that some $f \in
 \Gamma$ has positive entropy with respect
to some invariant  measure. Then $f$ must have somewhere non trivial
  negative and positive Liapunov spaces $E^1$ and $E^2$. Observe
then that
these spaces must be isotropic (in the sense of the Lorentz metric)
 and hence
are 1-dimensional. At most only one of these directions, say $E^1$,  is
contained in ${\cal F}$. However, by the uniform attraction
of ${\cal F}$ (Proposition \ref{south.north}), the direction $E^2$
 is mapped by powers of $f$, near ${\cal F}$, which contradicts the fact
that it is preserved by $f$.

\paragraph{The case of partially hyperbolic groups.}

In order to estimate its boundary, it
is easy to see that  we may replace $\Gamma$
 by a subgroup $\{ g^n/ n \in {\bf Z} \}$
generated by an element $g$ as in  Theorem \ref{theo.hyperbolic}.
We will thus prove:
$\partial_\infty \Gamma = \{ {\cal L}^1, {\cal L}^2 \}$
(following the notations of Theorem \ref{theo.hyperbolic}). More
precisely,  we will prove that:
$\lim_{n \to +\infty}g^n = {\cal L}^2$ and
$\lim_{n \to -\infty}g^n = {\cal L}^1$.

Let $( f_n )$ be a sequence of $\Gamma$
of the form $f_n = g^{k_n}$, for
some sequence of integers $(k_n)$. Suppose that $(f_n)$
has an approximately stable foliation, that is,  $( f_n^{-1} )$
converges in $\Gamma \cup \partial_\infty \Gamma$ to some foliation
${\cal F}$.

Suppose that $k_n \to +\infty$, when $n \to +\infty$. Then,
 by localization to larger and larger compact subsets, and using Theorem
\ref{theo.hyperbolic}, one deduces that the limit of $( g^{k_n} )$ cannot
 be different from ${\cal L}^2$. Therefore,
by compactness, $g^n \to {\cal L}^1$, when $n \to +\infty$ (because
${\cal L}^2$ is the unique limit of its convergent subsequences).
 By the same argument $\lim_{n \to -\infty}g^n = {\cal L}^1$.

It then follows that any sequence $(f_n= g^{k_n} )$
 for $(k_n)$ oscillating between $-\infty $
and $+\infty$ is not convergent. Hence,  the boundary
of $\Gamma$
is $\{ {\cal L}^1, {\cal L}^2 \}$.

\paragraph{Proof of Theorems \protect\ref{theo.isometry} and
\protect\ref{theo.flow}.}

Let $f  $ be an isometry of $M$ generating a
non-equicontinuous subgroup $\Gamma$.   Firstly, observe that
$\Gamma$ is closed. Indeed,  the closure of $\Gamma$
 is an abelian Lie group and hence, up to finite index,  it can
be written as a product $ \overline { \Gamma} = {\bf T}^i\times {\bf R}^j \times
{\bf Z}^k$, where ${\bf T}$  is  the torus part. But,
because $\overline{\Gamma}$ is non compact and  has  ${\bf Z}$
as a dense subgroup, then we must have $\overline{\Gamma} = {\bf Z}$.
Therefore,  $\Gamma$ is closed and amenable, and
hence elementary.

We apply  Theorem \ref{theo.elementary} and we get,
with the help of the previous notations,
 ${\cal L}^1$ (resp. ${\cal L}^2$).
 as an approximately
stable (resp. unstable) foliation  for $f$.

The  weak partial ergodicity
part of
 Theorem \ref{theo.isometry} follows from
the analogous statement in Theorem
\ref{theo.generalized}.
This ends the proof of Theorem \ref{theo.isometry}.

The same argument yields approximately stable
and unstable foliations for
  isometric flows, as stated in Theorem \ref{theo.flow}. The statement
 concerning the causal character of
non-equicontinuous flows, that is, their infinitesimal generators
must be non-timelike, was noticed in \cite{Zeg.Killing}.



\section{Non bi-polarized manifolds. Proof of Theorem \protect\ref{theo.tas}}


Observe that, over ${\cal U}$,
 $TAS$ is Lorentzian, that is the restriction
of the  metric along
$TAS$ is  of Lorentzian type.  Indeed, for
$x \in {\cal U}$, $TAS_x$
contains at least two different isotropic
directions.

Let $K$ be a compact subset of ${\cal U}$
 over which $E$
is  continuous.  We will consider
as in \ref{def.non.escaping}, the non escaping
set subset $NE(K, (f_n))$,
where $(f_n)$ is sequence in $\hbox{\rm Isom}(M)$.
 We firstly present the following
relative version of Corollary \ref{continuity.field}:

\begin{fact} \label{continuity.field.relative}
 Let $K$ be a compact subset of ${\cal U}$
 over which
$E$
is  continuous.
 Moreover, let
  ${\cal F} = \lim f_n \in
\partial_\infty\hbox{Isom}(M)$, and $X$ (resp. $Y$)
a vector field tangent to
${\cal F}^\perp$ (resp. ${\cal F}$). Then,  there
is a sequence of continuous
sections  $(X_n)$ (resp. $(Y_n)$)
of
$E$ over $K$, such that $X_n  \to X$ (resp.
$Y_n \to Y$) in the $C^0$ topology, and
satisfying the
following. Choose
an auxiliary norm
$\vert. \vert$ on $TM$. Then, there is
a real sequence $(a_n)$
 converging to $0$ (resp. a
bounded real sequence $(b_n)$)
such that
$\vert Df_nX_n(x_n) \vert <a_n  $
(resp. $\vert Df_nY_n(x_n)\vert <b_n$
whenever $x_n \in K \cap f_n^{-1}K$.
\end{fact}

\begin{fact} Let ${\cal F} \in \partial_\infty \hbox{Isom}M$.
Then:
$t(u, v) = 0$ for
$u \in T_x{\cal F}^\perp$,
$v \in T_x{\cal F} \cap E_x$, and
 almost every
$x \in {\cal U}$,

\end{fact}

\begin{proof}
Let $K$ be a compact subset over which
$E$
 and $t$ are  continuous.
Let $x \in NE(K, (f_n))$ and $u$, $v $ as
in the Fact.
 Extend
$u$ and  $v$ to local vector fields $X$
and
$Y$ tangent
 to ${\cal F}^\perp$
 and ${\cal F}$,  respectively. Approximate
 $X$ and $Y$ as in  Corollary
\ref{continuity.field}, and consider a sequence $x_n \to
x$ in $K$, such that $f_nx_x \in K$.
Thus, by continuity,  $t(u, v)
= \lim_{n \to \infty}t(X(x_n), Y(x_n))$
which equals $0$, from  the properties of
$X_n$ and $Y_n$
in  Corollary \ref{continuity.field}.

Finally,  recall that (see \ref{def.non.escaping})
$\hbox{Vol}NE(K, (f_n)) >
 \hbox{Vol}({\cal U})-2\epsilon$,
whenever $\hbox{Vol}(K) >
\hbox{Vol}({\cal U})$. Therefore,
the property holds is almost everywhere
because we can choose the volume
of $K$ arbitrarily approaching that of
${\cal U}$.
\end{proof}

It then follows   that
$t(u, v) = 0$ if $u \in S_x$ and
$v \in u^\perp \cap E_x$, for almost
every $x \in {\cal U}$.

Let $\{e_i \}$ be a basis
of $ T_xM$, and
 write $t(X, Y) = \Sigma
b_i(X, Y)e_i$ where the $b_i$
are bilinear scalar forms. To prove  Theorem \ref{theo.tas},
it suffices to show
 that for  each $b_i$
there is $\alpha_i$ such that
$b_i(u, v) =\alpha_i<u, v>$, whenever $u
\in TAS_x$.

So let $b$ be one of this forms
and write it
 $b(u, v) = <Au, v>$
 for some linear endomorphism of $E_x$.
From the above, one sees
that every
$u \in S_x$, is an eigenvector for $A$, $Au = \lambda_uu$.
 Now,
if $\lambda_u$ does not depend on $u$,
 $A$ is
 a homothety on $TAS_x$ and we are done.
 If not,  $A $ induces an eigenspace decomposition
on $E_x$.

By considering all
 the endomorphisms corresponding to
the $b_i$, and by letting $x$ varying  over ${\cal U}$, we
get an invariant measurable  decomposition
$E = E^1 \bigoplus \ldots \bigoplus E^k$
 of $TAS$ such that
$S_x \subset  E^1_x \cup \ldots \cup E^k_x$.
To show that  this is
impossible (and hence finishing the proof
of the theorem) we  use
 the following irreducibility fact:

\begin{fact}  There are no $\hbox{Isom}(M)$-invariant
measurable subbundles $E^1, \ldots, E^k$
of $TAS$, with $k \geq 2$, such that
  $S_x \subset
E^1_x  \cup \ldots \cup E^k_x$ in
some subset of ${\cal U}$
with positive volume
(it
is here where we use that
$\hbox{Isom}(M)$ is non elementary).
 \end{fact}

\begin{proof} We may suppose $k=2$ and
that $E^1$ is Lorentzian in
some subset of positive volume.
Since  we are  dealing  with measurable bundles, we
may suppose this is
everywhere true, just  restricting
 domains of definition.
So, we have an orthogonal decomposition
$TM = E^1 \bigoplus {E^1}^\perp$.

 Let ${\cal F} = \lim f_n \in
\partial_\infty \hbox{Isom}(M)$ and apply
 Fact
\ref{continuity.field.relative}, to a compact $K$
over which
$E^1$ and $E^2$ are continuous.  Observe that
if $u_n \in {E^1}^\perp_{x_n} \to u \neq 0$, and
$f_nx_n \in K$, then  $(Df_nu_n)$
 can't tend to $0$ because the metric on
${E^1}^\perp$ is Riemannian (i.e. positive definite).
This implies as in the proof
of the Fact above, that, over
$NE(K, (f_n))$, ${\cal F}^\perp$ is contained in
$E^1$.
Choosing $K$ with larger volume as necessary,
  we may  conclude that this is always true.
 \end{proof}

To finish the proof of Theorem \ref{theo.tas}, it
 remains to check its last partial ergodicity
statement. Let $\sigma$ be a function as
in the theorem. From
the ergodicity
property in Theorem \ref{theo.generalized},
$\sigma$ is constant along
any (1-dimensional) foliation
defined by a vector field $X \in S$ (see
the notations of \S \ref{ergodic.tas}).  Therefore,
by definition
$d \sigma\vert TAS = 0$.

\section{Examples} \label{examples}

Here we will give some examples of boundaries of
isometry groups of compact Lorentz manifolds. As we will
see in the part II of this work, the non-trivial
(i.e. non bi-polarized)
 cases,
involve constant curvature manifolds.  In what
follows, we therefore
investigate
 the structure of isometry groups
of manifolds which are  locally isometric to
a product of a Riemannian manifold by a
 constant curvature Lorentz manifold.
The following theorem summarizes the non trivial cases.

\begin{theorem} \label{theo.example}
Let $M$ be a compact Lorentz
manifold whose universal cover is a
product of a simply connected
Riemannian manifold $\tilde{N}$ and a complete
simply connected space $X_c$ of constant curvature
$c$. Let $\pi  \subset \hbox{Isom}\tilde{N} \times
\hbox{Isom}(X_c)$. Then,

(i) Such a manifold does not exist if  $c>0$.

(ii) If $c>0$ and  $M$ is not bi-polarized, then
$dimM =3$, and $M$ is
a quotient of $\tilde{N} \times
\widetilde{SL(2, {\bf R})}$ by
a subgroup
of $\hbox{Isom}\tilde{N} \times
\widetilde{SL(2, {\bf R})}$. The isometry
group of
$M$ is
a product of a compact group
by a finite cover of $PSL(2, {\bf R})$.
 Its boundary is the circle
$S^1$
endowed with the usual action of this latter group.

(iii) If $c = 0$ and $M$ is not bi-polarized, then
there is a metric decomposition $\tilde{M}=
\tilde{N^\prime} \times {\bf R}^d$, where
$\tilde{N^\prime}$ is a  Riemannian manifold, and ${\bf R}^d $ is a
Minkowski space, and such that
$\pi = \pi_1(M) \subset \hbox{Isom}\tilde{N^\prime}
\times {\bf R}^d$. The action of $\hbox{Isom}(M)$
on its boundary factors trough the action of
a (arithmetic) lattice of a $(d-1)$-dimensional
hyperbolic space
on its  sphere at infinity.

\end{theorem}

\begin{remark} [Erratum]
{\em  The case (iii), i.e. the flat
non-bi-polarized case,
was forgotten in
\cite{Zeg.CRAS}. So, the statement of the principal result
in this reference has to be modified to take
into account this case, see the part II of the
present article.

}

\end{remark}

Observe that that  we don't  address  completeness
questions, which will
be treated in detail in the part II of this article.
The theorem  will be proved through the present  section, which
contains further details.
Let's start with the following general fact
about de Rham decomposition of Lorentz manifolds,
which follows from the foliated version \ref{theo.foliated.version}.

\begin{proposition} Let $M$ be a compact Lorentz
manifold whose universal cover admits
a de Rham decomposition $\tilde{M} =
 \tilde{N} \times X$, where
$\tilde{N}$
 is Riemannian and
$X$ is Lorentzian. Then any foliation ${\cal F}$,
element of the boundary of $\hbox{Isom}(M)$
lifts to a foliation
$\tilde{ {\cal F}}$ containing the factor $\tilde{N}$, that
is its leaves have the form
$\tilde{N} \times \tilde{{\cal L}}_x$,
where $\tilde{{\cal L}}$
is a codimension one lightlike
geodesic foliation
of $X$, invariant under the action
of the projection of
$\pi_1(M)$
in
$\hbox{Isom}(X)$.

\end{proposition}

\subsection{De Sitter   space}
Let ${\bf R}^{p,q}$ be the space ${\bf R}^{p+q}$
endowed with a non degenerate quadratic form
of signature
$(-p, q)$ (for  example ${\bf R}^{1,q}$ is endowed
with a form of signature $-+\ldots+$). For
a real $r$, we denote by ${\bf S}^{p, q}(r)$ the
level $r$ in ${\bf R}^{p, q}$. The (universal)
de  Sitter space of dimension $q$
is ${\bf S}^{1, q}(+1)$. We only consider
the case
$q>2$, and hence, the de Sitter space  is
simply connected. It inherits from
${\bf R}^{1, q}$ a Lorentz metric
of positive constant curvature, and has isometry
group
$O(1, q)$.

The  well known Calabi-Markus phenomena
states
that a Lorentz manifold covered by
 ${\bf S}^{1,q}(+1)$
 has  a finite fundamental group. Let
 $\tilde{N}$
 be a complete simply connected manifold
and consider
the product $\tilde{M} = \tilde{N} \times
{\bf S}^{1,q}(+1)$.
This may have quotients with large fundamental
group, but
no compact ones. The proof of this claim
resembles
that of the Calabi-Markus phenomena and
 goes as follows. Suppose that $\pi \subset
 \hbox{Isom} \tilde{N} \times O(1, q)$ is
the fundamental
 group
of  a compact manifold $M$. Then,  the
 projection of
$\pi$ on $\hbox{Isom}\tilde{N}$ is not
 discrete, otherwise,
$M$ would fiber over a quotient of
$\tilde{N}$, with fiber
a compact quotient  of ${\bf S}^{1,q}(+1)$,
which does not exist. Hence,  there is a
divergent sequence
$ \gamma_n = (g_n, h_n) \in \pi$ such
that
$g_n \to 1$ in $\hbox{Isom}\tilde{N}$. The basic
 fact behind
the Calabi-Markus phenomena is that if
$K = {\bf S}^{1,q}(+1) \cap {\bf R}^{0, q}$, then
 $h(K) \cap K \neq
\emptyset$ for any $h \in O(1, q)$ (because
both $K$ and $h(K)$ are traces
 in ${\bf S}^{1,q}(+1)$ of
 linear  hyperplanes in
${\bf R}^{1+q}$). Therefore,  for any
open set $U \subset \tilde{N}$,
$\gamma_n(U \times K)
\cap (U\times K) \neq \emptyset$,
because $g_n \to 1$. This means
that $\pi$ does not act properly
on
$\tilde{N} \times {\bf S}^{1,q}(+1)$.

\subsection{Anti de Sitter manifolds}
The anti de Sitter space of dimension $1+q$ corresponds
to the level $-1$ in ${\bf R}^{2, q}$. More precisely,  the
(universal) anti de Sitter space ${\bf H}^{1, q}$ is
the universal cover of ${\bf S}^{2,q}(-1)$. However,
it is more convenient to work in the ``linear
 model'' ${\bf S}^{2,q}(-1)$ and then,  translate into
${\bf H}^{1,q}$. In fact one would hope that,  one
needs to pass to ${\bf H}^{1, q}$
only in pathological situations. Indeed,  as we will recall
below compact quotients of ${\bf H}^{1,q}$ are
in fact quotients of
${\bf S}^{2,q}(-1)$. So let's work in this latter space.
Let $Q = -x_1^2-x_2^2+x_3^2+\ldots x_{q+2}^2$
 be a  quadratic form defining ${\bf R}^{2,q}$.
Then ${\bf S}^{2,q}(-1)$ inherits a Lorentz metric
of negative constant curvature and its
isometry group
 is the orthogonal group of
$Q$, that is $O(2,q)$. The totally geodesic subspaces of
${\bf S}^{2,q}(-1)$ are exactly the traces on it of
linear subspaces of ${\bf R}^{2+q}$.
Furthermore, the lightlike
geodesic hypersurfaces
of ${\bf S}^{2,q}(-1)$
have the form:
$H_u= {\bf R}u^\perp \cap {\bf S}^{2,q}(-1)$
for $u \in {\bf R}^{2, q}$
an  isotropic vector.

Two hypersurfaces $H_u$
 and $H_v$ (for $u$ and $v$ isotropic)
 are disjoint if and only if
$u$ and $v$ are orthogonal but not
collinear. It follows that a codimension
one
lightlike geodesic foliation ${\cal F}$
 of ${\bf S}^{2,q}(-1)$ is determined
by hypersurfaces
$H_u$, for $u$ running over an
isotropic 2-plane $P$. The group
$O(2, q)$ acts transitively
on the space of isotropic 2-planes.
Its action on the space of pairs $(P_1, P_2)$
of isotropic 2-planes has exactly 3
orbits, one for $P_1 = P_2$, one
for $\hbox{dim}(P_1 \cap P_2) = 1$,
and the last
for $P_1 \cap P_2 = 0$.

We can write ${\bf R}^{2, q} = {\bf R}^{2, 2}
\bigoplus {\bf R}^{0, q}$. Because
the factor ${\bf R}^{0, 2}$ is Riemannian, most
of the dynamics happens  on ${\bf R}^{2, 2}$. Note
 that the decomposition above is not canonical,
but we will see this does not matter.
 So, we  now decorticate the case $q=2$, i.e.
the 3-dimensional anti de Sitter space. Instead of the standard
form $Q$,  one consider the form
$Q^\prime$ on ${\bf R}^4 = {\bf R}^2 \times {\bf R}^2$
defined by $Q^\prime(u, v) =
\omega(u, v)$, where $\omega$
is the volume form on ${\bf R}^2$. An element
$A \in SL(2, {\bf R})$
acts diagonally: $A(u, v) = (Au, Av)$,
by preserving  $Q^\prime$.
The group
$SL(2, {\bf R})$ acts freely on the
non-vanishing levels of
$Q^\prime$. Hence,  the $SL(2, {\bf R})$
orbits coincide (for dimension reasons)
with the components of  non vanishing levels. It turns out that
 metrics  defined in that way on
$SL(2, {\bf R})$
 are multiples of its Killing form. Indeed,  one
verifies that the identity component of
$O(2, 2)$ contains another copy of $SL(2, {\bf R})$
commuting with the given $SL(2, {\bf R})$-action. More
precisely, this identity component is
 isomorphic to $SL(2, {\bf R}) \times
Sl(2, {\bf R})$, which turns out to be
the identity component of the isometry group
of the Killing form
of $SL(2, {\bf R})$,
 acting by: $(g,h)x= gxh^{-1}$, where $g$, $h$, and $x$ belong
to
$SL(2, {\bf R})$.

 A lightlike
geodesic foliation
is determined by the orbits of a subgroup
 of the form $A\times \{1\}$
or $\{1\} \times A$ where $A$ is conjugate to
the affine group $AG \subset SL(2, {\bf R})$.
The case of
two different subgroups lying in the same factor
$\{1\} \times SL(2, {\bf R})$ or
$SL(2, {\bf R}) \times \{1\}$ corresponds to two
isotropic 2-planes such that $P_1 \cap P_2 =  0$.
The
case of different factors corresponds to  two
intersecting 2-planes.
It then follows that the
symmetry group
of a pair or isotropic 2-planes $(P_1, P_2)$
is conjugate up to swith of factors to either
$SL(2, {\bf R})\times\{1\}$ or $AG \times AG$.

For example,  the group $SL(2, {\bf R})$
acting as above on
$({\bf R}^4, Q^\prime)$ preserves each isotropic
2-plane of the form $\{(u, \alpha u), u \in {\bf R}^2 \}$
for $\alpha$ a real number, together with the plane
$\{0\} \times {\bf R}^2$. The set of the so defined
 foliations
is a circle, on which the other factor $SL(2, {\bf R})$
acts as usually.

From this we deduce  that the symmetry group
of 3 distinguish isotropic 2-planes is
contained up to swith of factors in
$SL(2, {\bf R})\times \{1\}$. In particular,
this symmetry group
is centralized by the other factor $\{1\} \times
 SL(2, {\bf R})$. Therefore, if a compact quotient $M$
of a product $\tilde{N} \times {\bf S}^{2,q}(-1)$,
is not bi-polarized and has non-compact
isometry group, then $\hbox{Isom}(M)$
contains $SL(2, {\bf R})$. From \cite{Gro}
(see also \cite{A-S.1} and \cite{Zeg.espacetemps.homogenes}),
$\hbox{Isom}(M)$
is a product of a compact group by a finite cover of
$PSL(2, {\bf R})$.

In higher dimensions the symmetry group of
3 isotropic 2-planes is contained in a product of
$SL(2, {\bf R}) $ by a compact group. This does not
act co-compactly on ${\bf S}^{2,q}(-1)$. By
translating the argumentation to
the universal cover ${\bf H}^{1, q}$, one
 checks that a co-compact subgroup cannot
preserve
3 distinct codimesion one lighlike geodesic
 foliations.
This proves the claim that if a compact Lorentz
 manifold
has a non-compact isometry group
and is not
bi-polarized, then the anti de Sitter factor of
its universal cover has
dimension 3.

\paragraph{Finitness of levels.}
As mentionned above, a compact  Lorentz manifold $M$, which
is a quotient  of
the universal anti de Sitter space,  is in fact up
to finite covers,
a quotient of the more concrete one
${\bf S}^{2, q}(-1)$.  This property was called
in \cite{K-R}
  the finitness of levels of compact anti de Sitter
 manifolds.
It is related to isometry groups as follows.
The fundamental group
of ${\bf S}^{2, q}(-1)$ is cyclic, generated by an element
$\sigma$.
The statement is that some power of $\sigma$
belongs to $\pi_1(M)$. Because
$\sigma$
is central in  $\hbox{Isom}({\bf H}^{1, q})$, it
defines an isometry $f$ of $M$. The finitness of the level
of $M$ is equivalent to
$f$ having finite order. Thus the opposite situation
 implies $\hbox{Isom}(M)$
 is infinite. If the  identity component of
$\hbox{Isom}(M)$
is not trivial, we get a connected subgroup
of $\hbox{Isom}({\bf H}^{1, q})$ centralizing
$\pi_1(M)$.
 One can thus  completely understand
this latter group, and in particular get a contradiction
to the hypothesis that  $f$ is of infnite order.
If $\hbox{Isom}(M)$ is discrete, then it is non-compact,
and therefore, it preserves a foliation, and thus
$\pi_1(M)$ is contained in the symmetry group of
a foliation as discribed above. Again here, by
working algebraically, one
gets a contradiction to the fact that $f$
has an infinite order.

\subsection{The flat case}
Let $g$ be a Lorentz quadratic form on
${\bf R}^k$.  It is called
rationnal, if some multiple $\alpha g$ has
rational coefficients when expressed
in the canonical basis. The isometry group
of the Lorentz space $({\bf R}^k, g)$
is the semi-direct product ${\bf R}^k
\rtimes O(g)$ where $O(g) \subset
GL(n, {\bf R})$ is the orthogonal group of $g$.
 Of course all these spaces are isometric to
the Minkowski space. However, this representation
may help to understand the modulus of flat
Lorentz structures.  For instance,
consider
  a topological torus ${\cal T}^k
 = {\bf R}^k/{\bf Z}^k$, then the familly
$({\bf T}^k, g)$, for $g$ as above,
exhautes the space of Lorentz flat
structures on ${\bf T}^k$.

A lightlike geodesic foliation in ${\bf R}^k$
consists obviously of  parallel lighltlike
affine hyperplanes, and  hence it is
determined by specifying  an isotropic direction
of $g$. Let $C(g)$ be  the space of such directions. It is
a sphere of dimension $k-2$. The action
of ${\bf R}^k \rtimes  O(g)$ on it,  factors via
the usual conformal action of $O(g)$ (i.e.
it
is isomorphic to the conformal action of $O(1, k-1)$ on
$S^{k-2}$).

Let  $M$ be a flat  Lorentz manifold
obtained as a quotient of $({\bf R}^k, g) $
 by a subgroup $\pi \subset {\bf R}^k \rtimes  O(g)$.
The  space $ {\cal  F G}(M)$ of codimension 1 lightlike geodesic
foliations of $M$, is identified
  with  the fixed points of the action of
$\pi$ (via its linear part) on $C(g)$.
It  follows in particular, that in the torus case,
this
space is identified with $C(g)$, for any $g$.
 Conversely ,  if $ \hbox{card}( {\cal F G} (M)) \geq 3$,
then  up to a subgroup
of finite index, $\pi$ is a lattice in the
translation part ${\bf R}^k$
of ${\bf R}^k \rtimes  O(g)$.  To see this, let
$E$ be  the linear  space generated by the
isotopic directions determined  by the elements
of ${\cal F G}(M)$.  Because of
the
dimension,
the quadratic form $g$
is of Lorentz type on
$E$ and is positive
 on $E^\perp$.

Observe now that
if an element
$A \in O(g)$
fixes $3$ isotropic directions $X_1$,
$X_2$ and $X_3$ , then  $A=
\pm \;  \hbox{Identity}$
on
on the linear space $F$
that they generate.
Indeed, write $AX_i = \lambda_iX_i$, then
$\lambda_i \lambda_j =1$ for
$i\neq j$, because $<X_i,X_j> \neq 0$,
 and hence $\lambda_i = \pm 1$.

It then follows that $\Gamma$
preserves a positive scalar
product. Therefore, by the
 Bieberbach Theorem,
$\pi$ is virtualy a lattice in
${\bf R}^k$. Equivalently
$M$ is covered by a torus.

As a corollary, we  get
that a compact flat manifold is
bi-polarized, unless it is, up to a	 finite cover,  a torus.
Now for a
torus $({\bf T}^k, g)$, its
isometry group is
$O(g, {\bf Z}) \rtimes {\bf R}^k$,
where $O(g, {\bf Z}) = O(g) \cap
GL(k, {\bf Z})$.

The action
of
$\hbox{Isom}({\bf T}^k, g)$
 on its boundary  can be  identified with
that of $O(g, {\bf Z})$
 on its limit set in
$C(g)$. So, we are
not leaving Kleinian  groups, but what
kind of group
could be $O(g, {\bf Z})$
and what it might have as  a limit set?

A classical
theorem of Harish-Chandra and Borel
states that $O(g, {\bf Z})$
 is  a lattice in
$O(g)$ if $g$ is rationnal. Hence its
limit
set in this case is the whole of
$C(g)$.
Let's now treat the general case (that is
when $g$ may be partially rational).

\begin{fact} Let $G$ be a non-compact  connected Lie  subgroup of
$O(g)$. Then either
$G$ fixes a point of $C(g)$ (that is $G$ is elementary) or
$G$ is reductive. In this last case, there is a Lorentzian
plane $E $ in ${\bf R}^k$ (that is $g$ restricted to
$E$ is of Lorentzian type) which is preserved by $G$
and such that $G = O(g \vert E) \times K$, where
$K$ is a compact subgroup
of $O(g \vert E^\perp)$. (Here for a non degenerate plane
$F$, we denote $O(g \vert F) =
\{ A \in O(g)/ A(F) = F \;  \hbox{and} \; A \vert F^\perp =
 \hbox{Identity} \}$).

\end{fact}

\begin{proof} The proof is standard and   briefly the idea is
as follows.
The radical $R$ of $G$ fixes a point of $C(g)$. If it is non
compact, $R$ contains a parabolic or hyperbolic
one parameter group. It then
follows that $R$ fixes at most two points, and in particular,
$G$ also fixes these points, and is thus elementary. Therefore,
if $G$ is not elementary, then its radical is compact.  Hence,
it is by definition reductive. The remaining part of the Fact is standard.
\end{proof}

Let $G$ be  the Zariski closure
  of $O(g, {\bf Z})$. It
is defined over ${\bf Q}$ (see for
instance \cite{G-PS}).   Therefore,  if non elementary,
$O(g, {\bf Z})$
is a lattice in $G$, by the Harish-Chandra-Borel
 Theorem, for general reductive
${\bf Q}$-groups.  Thus  $O(g, {\bf Z})$
may be thought as a lattice of
$O(g \vert E)$, and so
 its limit set is the projective
isotopic cone of $g \vert E$.

\paragraph{Partial rationality.}
In fact,  it seems that we have more precise
partial rationality
when $O(g, {\bf Z})$ is non elementary, that
is $E$ and $E^\perp$ seem to be rational. We were able
to check it, assuming $O(g, {\bf Z})$
non-cocompact in $G$. Indeed in this case, $O(g, {\bf Z})$
possesses  parabolic elements, that is
there exist an element of $O(g, {\bf Z})$
 of the form  $A = (B, C) \in O(g \vert E) \times K$
such that $B $ is unipotent, i.e. $B-1$ is nilpotent. In
fact unipotent elements of orthogonal
groups of Lorentz forms, have
degree of nilpotency $3$. Let ${\cal A} $
be the set of  $A \in O(g, {\bf Z})$ such that $B$
is unipotent.

Then
$E$ is contained in the vector
subspace $F = \cap_{A \in {\cal A}} \ker (A-1)^3$,
which is rational.

Observe that if $A \in {\cal A}$, then its projection
$C$ in $K$, acts trivially on
$ F^\prime= F \cap E^\perp$, because there it
detemines an unipotent element in a compact group.
In other words,  $F^\perp \subset \cap_{A \in {\cal A}}
(\ker (A-1)$. In fact,  we have equality.
  Otherwise,  the intresection $E \cap (\cap_{A \in
{\cal A}} \ker(A-1)$
would give a proper subspace of $E$ invariant
by $O(g, {\bf Z})$; which is easily seen
to be impossible. It then follows that $F^\prime$
 is rational. By similar arguments,  one
sees that $E = \bigoplus_{A \in {\cal A}} (A-1)(F)$,
and is hence rational. Again, consider $
F^{\prime \prime} = \bigoplus_{A \in {\bf A}}
\hbox{Image}(A-1)^3$.
This is  a rational space, which
is (by similar argument)
a supplementary of $F^\prime$ in $E^\perp$.
Hence $E^\perp$
 is rational.

In conclusion, at least when
$O(g, {\bf Z})$
is non elementary and  non co-compact, then up to a finite cover,
$M$
 is a metric product of a Lorentz rational torus
with a Riemannian torus.

\paragraph{The product case.}
 Suppose now that $M$ is a quotient
of a product $\tilde{N} \times {\bf R}^k$, where
$\tilde{N}$ is Riemannian and ${\bf R}^k$
 is a Minkowski
space. As above, assuming $M$
is not bi-polarized, we construct
 another invariant metric
decomposition $\tilde{M}
= \tilde{N}^\prime \times {\bf R}^d$, where
the projection of $\pi_1(M)$
on the isometry group
of the Minkowski space ${\bf R}^d$
consists of translations, only. Here
${\bf R}^d$
 correponds to the space $E$ defined above
 and $\tilde{N}^\prime = \tilde{N} \times
 E^\perp$.
Now we can write $\hbox{Isom}(M) =
\hbox{Nor}(\pi_1(M)) / \pi_1(M)$ where
$\hbox{Nor}(\pi_1(M))$
is the normalizer
of $\pi_1(M)$
in $\hbox{Isom}(\tilde{N}^\prime) \times
 \hbox{Isom}({\bf R}^d)$.
Nevertheless, it is not abvious how to
exploit further
such a formula in an algebraic way, and
so we  argue
as follows. Observe   that
the translations along ${\bf R}^d$
 centralize $\pi_1(M)$
and thus determine an isometric action
of ${\bf R}^d$
on $M$. This action is contained in a
compact group
because of its abvious equicontinuity. Therefore,
we get a torus ${\bf T}^k$ in
$\hbox{Isom}(M)$. This torus inherits  a
  metric defined on
its Lie algebra ${\cal T}^k$,  by means
of the formula:
$<X,Y> = \int <X(x), Y(x)>$,  where
$X$ and $Y$
are Killing fields generating
flows in ${\bf T}^k$.
 Parallel fields tangent to
${\bf R}^d$, allow us to embed
isometrically the Minkowski space
${\bf R}^d$
in $({\cal T}^k, <,>)$. In fact, as $\tilde{M}$
itself, ${\cal T}^k$
admits an orthogonal  decomposition
$A \bigoplus {\bf R}^d $, such that
 the scalar product $<,>$ on $A$
is positive definite, and hence the metric on
${\bf T}^k$
is Lorentzian. This construction is natural
in all its steps, and therefore we have succeded to
essentially incorporate questions about
the isometry group of $M$  into
ones relative to
${\bf T}^k$. In particular one sees
that
the action of $\hbox{Isom}(M)$
on its boundary factors through
the action of an arithmetic lattice
of the hyperbolic space ${\bf H}^{d-1}$
on its boundary ${\bf S}^{d-2}$.

\medskip
\noindent
CNRS, UMPA, \'Ecole Normale Sup\'erieure de Lyon \\
46, all\'ee d'Italie,
 69364 Lyon cedex 07,  FRANCE \\
Zeghib@umpa.ens-lyon.fr

\end{document}